\documentclass[preprint,12pt]{elsarticle}

\usepackage{amssymb}
\usepackage{amsmath}
\usepackage{lineno}

\usepackage{hyperref}
\usepackage{cleveref}
\usepackage{color,soul}
\usepackage{orcidlink}

\biboptions{sort&compress}

\newtheorem{thm}{Theorem}
\newtheorem{lem}[thm]{Lemma}
\newdefinition{rmk}{Remark}
\newdefinition{dfn}{Definition}
\newproof{pf}{Proof}
\newproof{pot}{Proof of Theorem \ref{thm2}}

\usepackage{tikz}
\usepackage{subfig}
\crefname{equation}{}{}

\begin{document}

\begin{frontmatter}

\title{Neural network methods for Neumann series problems of Perron-Frobenius operators}        

\author[1]{Tanakorn Udomworarat\orcidlink{0009-0002-9879-1885}\corref{cor1}}
\ead{Tanakorn.Udomworarat@nottingham.ac.uk}
\author[1]{Ignacio Brevis\orcidlink{0000-0003-1620-019X}}
\ead{Ignacio.Brevis1@nottingham.ac.uk}
\author[1]{Martin Richter\orcidlink{0000-0001-7958-7550}}
\ead{Martin.Richter@nottingham.ac.uk}
\author[2]{Sergio Rojas\orcidlink{0000-0001-7203-7740}}
\ead{Sergio.rojas@monash.edu}
\author[1]{Kristoffer G. van der Zee\orcidlink{0000-0002-6830-8031}}
\ead{KG.vanderZee@nottingham.ac.uk}

\cortext[cor1]{Corresponding author}

\affiliation[1]{organization={School of Mathematical Sciences},
addressline={University of Nottingham},
country={United Kingdom}}

\affiliation[2]{organization={School of Mathematics},
addressline={Monash University},
country={Australia}}

%% Abstract
\begin{abstract}
Problems related to Perron-Frobenius operators (or transfer operators) have been extensively studied and applied across various fields. In this work, we propose neural network methods for approximating solutions to problems involving these operators. Specifically, we focus on computing the Neumann series of non-expansive Perron-Frobenius operators under a given $L^p$-norm with a constant damping parameter in $(0,1)$. We use PINNs and RVPINNs to approximate solutions in their strong and variational forms, respectively. We provide a priori error estimates for quasi-minimizers of the associated loss functions. We present some numerical results for 1D and 2D examples to show the performance of our methods. We also demonstrate the applicability of our methods by approximating interior densities in a two-cavity system.
\end{abstract}

%% Keywords

\end{frontmatter}


\section{Introduction}
Perron-Frobenius operators, also known as transfer operators, play an important role in capturing behaviors of dynamical systems. These operators have been applied to various fields of research, including engineering \cite{vaidya2010nonlinear,hartmann2019high}, earth sciences \cite{froyland2014well}, and atmospheric sciences \cite{tantet2015early,tantet2018crisis}. Another viewpoint to describe system behaviors is studying observables using Koopman operators, which are the adjoints of Perron-Frobenius operators. While Koopman operators track the evolution of functions on the state space, Perron-Frobenius operators describe statistical behavior through the evolution of densities \cite{klus2015numerical}. These operators provide linear representations of (nonlinear) dynamical systems, which are important for estimating the long-term system behavior and prediction of future states \cite{kaiser2020data}.

One particularly interesting problem associated with Perron-Frobenius operators is finding their Neumann series when applied to a given initial density. The solution to this problem corresponds to the accumulated density in the long-time limit, converging toward the stationary density. This problem is closely related to the dynamical energy analysis approach \cite{tanner2009dynamical} used to determine wave energy distributions in the high-frequency regime, which has a direct application in finding equilibrium energy distributions \cite{slipantschuk2020transfer, richter2020convergence}. A classical approach for approximating solutions to problems involving such operators is Ulam's method \cite{ulam1960acollection}, which approximates Perron-Frobenius operators by projecting them onto a subspace spanned by characteristic functions (a space of piecewise constant functions). However, without additional effort, fixed-grid-based methods, including Ulam's method, are known to have limitations, particularly in handling irregularities, and the convergence rate is relatively low \cite{bose2000exact}. 

In recent years, neural networks have been applied more and more to improve the approximate solution of various mathematical problems (see \cite{yu2018deep,sirignano2018dgm,raissi2019physics,kharazmi2019variational,kharazmi2021hp,rojas2024robust,cai2021physics} for examples), including problems in dynamical systems (see \cite{faria2024data,linka2022bayesian,saqlain2023discovering,oluwasakin2023optimizing,antonelo2024physics,yu2024learning} for examples). Among these methods, Physics-Informed Neural Networks (PINNs) \cite{raissi2019physics} and various versions of Variational PINNs (VPINNs) \cite{kharazmi2019variational,kharazmi2021hp,rojas2024robust} have gained popularity for solving equations involving operators, such as partial differential equations (PDEs). PINNs were developed to generate neural network functions to approximate solutions by minimizing a loss function based on the PDE residual. Similar methods can be found in \cite{cai2020deep,brevis2022neural}. VPINNs \cite{kharazmi2019variational} were introduced, aiming to solve PDEs in variational forms, with further studies presented in \cite{berrone2022solving,berrone2022variational}. Robust VPINNs (RVPINNs) were developed to give a stable version of VPINNs by minimizing a loss function based on the discrete dual norm of the residual of the variational equation. Throughout this work, we use the terms PINNs and RVPINNs to refer to neural networks that encode equations describing systems by minimizing the residual of the equations in their strong and variational forms, respectively. 

Neural networks offer several advantages for approximating solutions, including high expressivity \cite{Gühring_Raslan_Kutyniok_2022}, handling complex domains \cite{cuomo2022scientific}, and scalability to high-dimensional problems \cite{onken2022neural}. Several papers have used neural networks to learn operators in dynamical systems. For example, the authors of the seminal work \cite{lusch2018deep} have proposed a deep learning framework to learn a finite-dimensional approximation of Koopman operators. However, few papers have been dedicated to studying neural networks for problems involving Perron-Frobenius operators. Applying neural networks to solve problems related to Perron-Frobenius operators, in particular Neumann series problems, would be beneficial for approximating irregular solutions or addressing high-dimensional spaces. 

In this article, we propose and analyze neural network methods for approximating solutions to Neumann series problems for Perron-Frobenius operators with an initial density in $L^p$, where $1<p<\infty$, considering both their strong and variational forms. Our primary focus is on Perron-Frobenius operators that are non-expansive under the $L^p$-norm with a constant damping parameter in $(0,1)$. For consistency and simplicity, we present most results for solutions in $L^2$ but also provide details for solutions in $L^p$. 

The main contributions of this work are as follows. We define an abstract setting for non-expansive Perron-Frobenius operators from $L^p$ to $L^p$. We establish the well-posedness of the variational formulation in $L^2$ and $L^p$ settings. We propose neural network methods for solving Neumann series problems based on PINNs and RVPINNs. A key advantage of the RVPINNs approach is that it does not require the inverse map of the underlying dynamical system, as highlighted in Remark \ref{remark: reformulated RVPINNs loss function}. To derive a priori error estimates for the proposed methods, we recall the concepts of neural network function manifolds and quasi-minimizers presented in \cite{brevis2022neural}. The local Fortin's condition modified from \cite{rojas2024robust} is required to analyze the RVPINNs method. We implement the proposed methods for approximating the Neumann series of Perron-Frobenius operators associated with the tent map, the boundary map in a circular domain, and the standard map to show the performance of our methods. We also demonstrate the applicability of our methods by approximating interior densities in a two-cavity system. We compare our methods against a fixed-grid-based approach and the truncated sum of Neumann series.

The article is organized as follows. In section \ref{sec: Abstract framework}, we introduce Perron-Frobenius operators, state our assumptions, present Neumann series problems in both strong and variational forms, and provide the well-posedness of the variational problems when solutions are in $L^2$. We also review fixed-grid-based methods, including Ulam's method, and provide an associated error estimate. Section \ref{sec: Neural network framework} describes neural network frameworks, including our proposed methods for approximating solutions, focusing on PINNs and RVPINNs. Section \ref{sec: Analysis of methods} is devoted to the analysis of these methods. Section \ref{sec: Numerical examples} presents numerical examples. Section \ref{sec: conclusions} is the conclusions and possible extensions. \ref{section: well-posedness} provides the well-posedness of the variational problem when solutions are in $L^p$. \ref{section: equivalent form of RVPINNs loss} shows the equivalent form of the RVPINNs loss function. \ref{section: error estimates proof} shows error estimate proofs for PINNs and RVPINNs.

\section{Abstract framework}\label{sec: Abstract framework}

Let $\Omega$ be a bounded subset of $\mathbb{R}^d, d \geq 1$, and $S:\Omega \to \Omega$ be a map describing a discrete dynamical system. Given a measure space $(\Omega,\mathcal{B},\mu)$ where $\mathcal{B}$ is a Borel sigma-algebra and $\mu$ is the Lebesgue measure, the Perron-Frobenius operator associated with $S$ is defined as the linear operator $\mathcal{P}: L^1(\Omega) \to L^1(\Omega)$ such that
\begin{equation}\label{eq: Perron-Frobenius definition}
    \int_A \mathcal{P}f d\mu = \int_{S^{-1}(A)}f d\mu, \quad \forall A \in \mathcal{B}, \forall f \in L^1(\Omega).
\end{equation}
In the special case where $\Omega=[a,b]$ is a closed interval, the Perron-Frobenius operator $\mathcal{P}$ can be expressed explicitly (see, e.g., \cite[equation 3.2.6]{lasota2013chaos}) as:
\begin{equation}\label{eq: Perron-Frobenius closed interval}
    \mathcal{P}f(x)=\dfrac{d}{dx}\int_{S^{-1}([a,x])}f(s)ds.
\end{equation}

It is well-known that the Perron-Frobenius operator $\mathcal{P}$ is non-expansive under the $L^1$-norm, i.e., $\|\mathcal{P}f\|_{L^1} \leq \|f\|_{L^1}$ for all $f\in L^1(\Omega)$. For simplicity\footnote{$L^p(\Omega)$ (with $1<p<\infty$) has convenient properties for analysis such as reflexivity and strict convexity, but these properties do not hold in $L^1(\Omega)$  \cite{muga2020discretization}.}, we consider densities $f$ in $L^p(\Omega)$, where $1 < p < \infty$, which is a subspace of $L^1(\Omega)$. So, we need additional assumptions on the Perron-Frobenius operator for defining our problem in $L^p$-setting.
In this article, we make the following assumptions:
\begin{enumerate}
    \item[(A1)] Density functions $f$ are nonnegative functions in $U:=L^p(\Omega)$ with $1<p<\infty$,
    \item[(A2)] $\mathcal{P}$ is non-expansive under the $U-$norm, i.e., for all $f \in U$, 
    \begin{equation*}
        \|\mathcal{P}f\|_{U} \leq \|f\|_{U}.
    \end{equation*}
\end{enumerate}

\begin{rmk}
    The assumption (A2) is equivalent to stating that $\mathcal{P}:U\to U$ is a bounded operator with operator norm $\|\mathcal{P}\| \leq 1$.
\end{rmk}

\begin{rmk}\label{Remark: S diffeomorphism}
    The assumption (A2) holds in many dynamical systems. For instance, if $S$ is an invertible nonsingular transformation, the Perron-Frobenius operator $\mathcal{P}$ associated with $S$ is given by
    \begin{displaymath}
        \mathcal{P}f(x) = f(S^{-1}(x))\left|J_{S^{-1}}\right|, \quad \textrm{ for }x \in \Omega,
    \end{displaymath}
    where $\left|J_{S^{-1}}\right|$ denotes the Jacobian determinant of the inverse map $S^{-1}$ (cf. \cite[Corollary 3.2.1.]{lasota2013chaos}).
    Therefore, if $\left|J_{S^{-1}}\right| \leq 1$, $\mathcal{P}$ satisfies the assumption (A2). Indeed, by applying the change of variable, we obtain
    \begin{align*}
        \|\mathcal{P}f\|_{U} 
        &= \left(\int_{\Omega} \left( f(S^{-1}) \right)^p \left| J_{S^{-1}} \right|^p d\mu \right)^{1/p} \\ 
        &\leq \left(\int_{\Omega} \left( f(S^{-1}) \right)^p \left| J_{S^{-1}} \right| d\mu \right)^{1/p}
        = \left( \int_{S^{-1}(\Omega)} f^p d\mu \right)^{1/p}
        = \|f\|_{U}.
    \end{align*}
\end{rmk}

\subsection{Neumann series problem}

We are interested in approximating the solution of the following problem: Given a nonnegative function $f_0 \in U$, find $u \in U$ satisfying
\begin{equation}\label{Eq: steady state equation}
    u - \alpha\mathcal{P}u = f_0,
\end{equation}
where $0<\alpha<1$ is a damping parameter. This is equivalent to finding the Neumann series expressed as
\begin{equation}\label{Eq: power series P}
    u = f_0 + \alpha\mathcal{P}f_0 + \alpha^2\mathcal{P}^2f_0 + \cdots.
\end{equation}

\begin{rmk}
    In applications, the Neumann series (\ref{Eq: power series P}) represents the equilibrium energy distribution resulting from a continuous injection of the initial energy distribution $f_0$ into the system, with energy reduced by a factor $\alpha$ each time it reaches the boundary. In a more general setting, the damping parameter $\alpha$ can be replaced by a suitable weight function $\alpha:\Omega \to [0,\infty)$.  This would mean more realistic energy losses for individual trajectories.
\end{rmk}

\subsection{Variational problem}

The variational version of problem \cref{Eq: steady state equation} is the following:
\begin{equation}\label{steady state variational problem}
    \textrm{Find } u \in U \textrm{ such that } b(u,v) = l(v) \quad \forall v \in V,
\end{equation}
  where $V$ is the dual space of $U$, $b(\cdot,\cdot)$ and $l(\cdot)$ are a bilinear form on $U \times V$ and a linear functional on $V$, respectively, defined by 
  \begin{equation}\label{eq: def of b and l}
      b(u,v) := \int_{\Omega}(u-\alpha\mathcal{P}u) v d\mu \quad \textrm{ and } \quad l(v) := \int_{\Omega} f_0 v d\mu.
  \end{equation} 
  \begin{rmk}
      The bilinear form $b(\cdot,\cdot)$ defined in (\ref{eq: def of b and l}) can be written in terms of the associated Koopman operator\footnote{The Koopman operator associated with the map $S$ is defined by $\mathcal{K}v=v\circ S$.} $\mathcal{K}$, by applying the duality identity\footnote{Let $1<p,q<\infty$ with $\dfrac{1}{p}+\dfrac{1}{q}=1$. If $\mathcal{P}:L^p(\Omega)\to L^p(\Omega)$ is the Perron-Frobenius operator and $\mathcal{K}:L^q(\Omega)\to L^q(\Omega)$ is the associated Koopman operator, then $\int_{\Omega}(\mathcal{P}u)v d\mu = \int_{\Omega}u(\mathcal{K}v) d\mu$.}, as follows:
      \begin{equation}\label{bilinear K}
          b(u,v) := \int_{\Omega}u (v-\alpha \mathcal{K}v) d\mu.
      \end{equation}
      This form will be used to reformulate the RVPINNs loss function in Remark \ref{remark: reformulated RVPINNs loss function}.
  \end{rmk}
  
  When $U = V = L^2(\Omega)$, the variational problem \cref{steady state variational problem} satisfies the conditions of the Lax-Milgram Theorem, ensuring the existence and uniqueness of the solution. The details regarding the well-posedness for $U = L^p(\Omega)$ and $V = L^q(\Omega)$, where $1 < p,q < \infty$ with $\dfrac{1}{p}+\dfrac{1}{q} = 1$, are provided in \ref{section: well-posedness}.
  
\begin{thm}[Well-posedness of variational problem in $L^2$]
  Given $U = V = L^2(\Omega)$, the variational problem \cref{steady state variational problem} satisfies the following.
  \begin{enumerate}
      \item Boundedness of $b$: $|b(u,v)| \leq (1+\alpha)\|u\|_{L^2}\|v\|_{L^2}, \quad \forall u,v \in L^2(\Omega)$,
      \item Boundedness of $l$: $|l(v)| \leq \|f_0\|_{L^2}\|v\|_{L^2}, \quad \forall v \in L^2(\Omega)$,
      \item Coercivity: $b(u,u) \geq (1-\alpha) \|u\|_{L^2}^2, \quad \forall u \in L^2(\Omega)$.
  \end{enumerate}
  Hence, it has a unique solution.
\end{thm}

\begin{pf}
    1. Boundedness of $b$: For any $u,v \in L^2(\Omega)$, by using the Cauchy-Schwarz inequality, the triangle inequality, and the assumption (A2), we have
    \begin{align*}
          |b(u,v)| &= \left|\int_{\Omega}(u-\alpha\mathcal{P}u) v d\mu \right| \leq \int_{\Omega}|(u-\alpha\mathcal{P}u) v| d\mu \leq \|u-\alpha\mathcal{P}u\|_{L^2}\|v\|_{L^2} \\
          &\leq \left(\|u\|_{L^2} + \alpha\|\mathcal{P}u\|_{L^2} \right)\|v\|_{L^2} \leq (1+\alpha)\|u\|_{L^2}\|v\|_{L^2}.
      \end{align*}
      
    2. Boundedness of $l$: Since $f_0,v \in L^2(\Omega)$, by using the Cauchy-Schwarz inequality, 
    \begin{align*}
        |l(v)| = \left|\int_{\Omega} f_0 v d\mu\right| \leq \|f_0\|_{L^2} \|v\|_{L^2}.
    \end{align*}
    
    3. Coercivity: For any $u \in L^2(\Omega)$, by using the Cauchy-Schwarz inequality and the assumption (A2), we have
    \begin{align*}
        b(u,u) &= \int_{\Omega}(u-\alpha\mathcal{P}u)ud\mu = \int_{\Omega}u^2d\mu - \alpha\int_{\Omega}(\mathcal{P}u)ud\mu \\
        &\geq \|u\|_{L^2}^2 - \alpha\|\mathcal{P}u\|_{L^2}\|u\|_{L^2} \geq \|u\|_{L^2}^2 - \alpha\|u\|_{L^2}^2 = (1-\alpha)\|u\|_{L^2}^2.
    \end{align*}
\end{pf}

\subsection{Fixed-grid-based method}
Let $U=V=L^2(\Omega)$. A classical way for solving problem \cref{steady state variational problem} is to discretize the problem based on Galerkin projections. To do that, we first consider a finite-dimensional subspace $V_M \subset V$ with a basis $\{ g_1,\ldots,g_M \}$. Next, we seek $u_M \in V_M$ such that 
\begin{equation}
    b(u_M,g_m) = l(g_m) \quad \forall m \in \{1,2,\ldots,M\}, \label{steady state variational problem finite}
\end{equation}
where $b(\cdot,\cdot)$ and $l(\cdot)$ are defined in \cref{eq: def of b and l}.
Then, as $u_M \in V_M$, we write $u_M=\sum\limits_{k=1}^M c_kg_k$ where $c_k$ is a coefficient associated to $g_k$. Finally, we solve the system of equations:
\begin{align*}
    \sum_{k=1}^M c_k \int_{\Omega}\left(g_k - \alpha(\mathcal{P}g_k)\right)g_m d\mu &= \int_{\Omega} f_0g_m d\mu, \quad m = 1,2,\ldots,M. 
\end{align*}
Solving this system is equivalent to solving the matrix equation $A\underline{c} = \underline{b}$, where the matrix entries of $A$ are given by $a_{mk} = \int_{\Omega}\left(g_k - \alpha(\mathcal{P}g_k)\right)g_m d\mu$ and the components of the vector $\underline{b}$ are $b_m = \int_{\Omega} f_0g_m d\mu$.

\begin{rmk}
    If the domain $\Omega$ is partitioned into $M$ subdomains $\omega_1,\ldots,\omega_M$ and $g_m$ is chosen to be the normalized characteristic function on a subdomain $\omega_m$, which forms an orthonormal basis of $V_M$, the fixed-grid-based method is equivalent to what so-called Ulam's method \cite{ulam1960acollection}. The matrix equation is equivalent to $(I_M - \alpha P_M)\underline{c} = \underline{b}$, where $I_M$ is the $M\times M$ identity matrix and $P_M$ is the matrix representation of $\mathcal{P}$ given by $P_M = [p_{mk}]_{M \times M}$ with
    \begin{equation*}
        p_{mk} = \dfrac{\mu(S^{-1}(\omega_m) \cap \omega_k)}{\mu(\omega_k)}.
    \end{equation*} 
\end{rmk}

The following theorem shows an error estimate for the fixed-grid-based method.

\begin{thm}
    Let $u_M \in V_M$ be the solution of problem \cref{steady state variational problem finite} and $u$ be the solution of problem \cref{steady state variational problem} for $U=V=L^2(\Omega)$. We have
    \begin{equation*}
        \|u-u_M\|_{L^2} \leq \dfrac{2}{1-\alpha} \inf_{v_m \in V_M} \|u-v_m\|_{L^2}.
    \end{equation*}
\end{thm}

\begin{pf} For any $v_m \in V_M$, by using coercivity, Galerkin orthogonality, and the boundedness of $b$, we have
    \begin{align*}
        \|v_m - u_M\|_{L^2}^2 &\leq \dfrac{1}{1-\alpha}b(v_m-u_M,v_m-u_M) \\
        &= \dfrac{1}{1-\alpha}\left[b(v_m-u,v_m-u_M) + b(u-u_M,v_m-u_M)\right] \\
        &= \dfrac{1}{1-\alpha}b(v_m-u,v_m-u_M) \\
        &\leq \dfrac{1+\alpha}{1-\alpha} \|v_m-u\|_{L^2}\|v_m-u_M\|_{L^2}.
    \end{align*}
    Dividing both sides by $\|v_m-u_M\|_{L^2}$ yields
    \begin{equation*}
        \|v_m - u_M\|_{L^2} \leq \dfrac{1+\alpha}{1-\alpha} \|v_m-u\|_{L^2}.
    \end{equation*}
    Applying the triangle inequality, we have
    \begin{align*}
        \|u-u_M\|_{L^2} &\leq \|u-v_m\|_{L^2} + \|v_m-u_M\|_{L^2} \\
        &\leq \|u-v_m\|_{L^2} + \dfrac{1+\alpha}{1-\alpha} \|v_m-u\|_{L^2} \\
        &= \dfrac{2}{1 - \alpha} \|u-v_m\|_{L^2}.
    \end{align*}
    Taking the infimum over all possible $v_m \in V_M$ yields the theorem.
\end{pf}

\begin{rmk}\label{remark: error bound piecewise linear approx}
    In the case where $V_M$ is the space of piecewise linear finite elements defined on a mesh of size $h$, the error bound for a solution $u$ in the Sobolev space $H^2$ is of order $O(h^2)$  (cf. \cite[Theorem 1.5]{ainsworth2000posteriori}).
\end{rmk}

\section{Neural network framework}\label{sec: Neural network framework}

To approximate the solution of \cref{Eq: steady state equation} or \cref{steady state variational problem}, we use a neural network $u_{\theta}:\mathbb{R}^d \to \mathbb{R}$ where $\theta \in \mathbb{R}^s$ represents the trainable parameters, including the networks' weights and biases. The $L$-layer neural network consists of an input layer, $L-1$ hidden layers, and an output layer. Each hidden layer is a composition of an affine transformation and a nonlinear activation function. Specifically, the $j^{th}$ hidden layer, for $j = 1,\ldots, L-1$, is defined as 
\begin{equation*}
    \underline{x}^{(j)} = \phi^{(j)}(\underline{x}^{(j-1)}) = \sigma(W^{(j)}\underline{x}^{(j-1)} + \underline{b}^{(j)}),
\end{equation*}
where $\sigma$ is an activation function, $W^{(j)}$ and $\underline{b}^{(j)}$ are weight matrices and bias vectors, respectively. The number of columns in $W^{(j)}$ corresponds to the number of neurons in the $j^{th}$ layer. The output layer follows the same structure as the hidden layers but omits the activation function, i.e.,
\begin{equation*}
    \phi^{(L)}(\underline{x}^{(L-1)}) = W^{(L)}\underline{x}^{(L-1)} + \underline{b}^{(L)}.
\end{equation*}
Thus, the $L$-layer neural network is expressed as
\begin{equation*}
    u_{\theta}(\underline{x}^{(0)}) = \phi^{(L)} \circ \phi^{(L-1)} \circ \ldots \circ \phi^{(1)}(\underline{x}^{(0)}),
\end{equation*}
where $\underline{x}^{(0)}$ represents the input to the neural network.
We denote $\mathcal{M}_n$ as the set of neural network functions with a total of $n$ neurons, based on the architecture described above:
\begin{equation}\label{Eq: Neural network manifold}
    \mathcal{M}_n := \{u_{\theta}:\mathbb{R}^d \to \mathbb{R}~|~\theta \in \mathbb{R}^s\}.
\end{equation}
Assume that $\emptyset \neq \mathcal{M}_n \subset U$. We aim to find $u_{\theta} \in \mathcal{M}_n$ such that $u_{\theta} \approx u$ where $u$ is the solution of \cref{Eq: steady state equation} or \cref{steady state variational problem}. To do that, we train the neural network to minimize a loss function using some optimization methods, such as the BFGS method or Adam \cite{kingma2014adam}.

\subsection{PINNs}
To approximate the solution of problem \cref{Eq: steady state equation}, we follow the PINNs method by using a loss function defined in terms of the residual of the equation. The loss function\footnote{This loss function is slightly different from the original PINNs as we do not take the square of the residual norm.} is given by:
\begin{equation}
    \mathcal{L}_{\mathrm{PINNs}}(u_{\theta}) := \| f_0 - u_{\theta} + \alpha\mathcal{P}u_{\theta} \|_{U}.\label{Eq: loss function PINNs}
\end{equation}
The approximate solution of \cref{Eq: steady state equation} is obtained by solving the optimization problem:
\begin{equation}
    \textrm{Find } u_{\theta^*} \in \mathcal{M}_n \textrm{ such that } \theta^* = \operatorname*{argmin}_{\theta \in \mathbb{R}^s} \mathcal{L}_{\mathrm{PINNs}}(u_{\theta}).\label{PINNs problems}
\end{equation}
\subsection{RVPINNs}\label{section: RVPINNs}
For a given discrete space $V_M \subset V$, we consider the following Petrov-Galerkin type discretization of \cref{steady state variational problem}: 
\begin{equation}\label{steady state variational problem discrete}
    \textrm{Find } u_{\theta} \in \mathcal{M}_n \textrm{ such that } b(u_{\theta},v_M) = l(v_M) \quad \forall v_M \in V_M.
\end{equation}
We follow the RVPINNs method to approximate the solution of problem \cref{steady state variational problem discrete} by defining a loss function based on the discrete dual norm of the variational equation residual. The loss function\footnote{This loss function is slightly different from the original RVPINNs as we do not take the square of the residual norm.} is given by:
\begin{equation}
    \mathcal{L}_{\mathrm{RVPINNs}}(u_{\theta}) := \sup_{0 \neq v_M \in V_M} \frac{l(v_M)-b(u_{\theta},v_M)}{\|v_M\|_V}.\label{Eq: loss function VPINNs sup}
\end{equation}
The approximate solution of \cref{steady state variational problem discrete} is obtained by solving the optimization problem:
\begin{equation}
    \textrm{Find } u_{\theta^*} \in \mathcal{M}_n \textrm{ such that } \theta^* = \operatorname*{argmin}_{\theta \in \mathbb{R}^s} \mathcal{L}_{\mathrm{RVPINNs}}(u_{\theta}).
\end{equation}
In practical settings, where $U=V=L^2(\Omega)$ and $V_M$ is a finite-dimensional subspace of $L^2(\Omega)$ with an orthonormal basis $\{g_m\}_{m=1}^M$, the RVPINNs loss function can be expressed in the following equivalent form:
\begin{equation}
    \mathcal{L}_{\mathrm{RVPINNs}}(u_{\theta}) =  \sqrt{\sum_{m=1}^M \left(\int_{\Omega} (f_0 - u_{\theta} + \alpha\mathcal{P}u_{\theta})g_m d\mu\right)^2}.\label{Eq: loss function VPINNs}
\end{equation}
The equivalence between these two forms is shown in \ref{section: equivalent form of RVPINNs loss}. 
It is worth noting that computing these loss functions requires numerical integration methods.

\begin{rmk}\label{remark: reformulated RVPINNs loss function}
    Using the bilinear form \cref{bilinear K}, we can reformulate the loss function \cref{Eq: loss function VPINNs} to eliminate the composition between the Perron-Frobenius operator $\mathcal{P}$ and the neural network $u_{\theta}$.  This reformulation avoids the need for the inverse map $S^{-1}$ associated with the operator $\mathcal{P}$, and instead only requires the forward map $S$ for evaluating the Koopman operator $\mathcal{K}$ applied to the predefined test functions $g_m$. The resulting RVPINNs loss function is given by:
\begin{equation*}
    \mathcal{L}_{\mathrm{RVPINNs}}(u_{\theta}) =  \sqrt{\sum_{m=1}^M \left(\int_{\Omega} f_0g_m - u_{\theta}(g_m - \alpha\mathcal{K}g_m) d\mu\right)^2}.
\end{equation*}

\end{rmk}

The following section presents a priori error estimates for PINNs and RVPINNs. 

\section{Analysis of methods}\label{sec: Analysis of methods}

It is known that $\mathcal{M}_n$ defined in \cref{Eq: Neural network manifold} is neither a linear, a closed, nor a convex subset of $U$ \cite{petersen2021topological}. Although there is an infimum of a loss function $\mathcal{L}$ in $U$, the minimizer in $\mathcal{M}_n$ may not exist \cite{brevis2022neural}. To study error analysis related to neural networks, we will use a relaxed definition of a minimizer called a quasi-minimizer (cf. \cite{brevis2022neural}).
\begin{dfn}
    Let $\mathcal{L}:U \to \mathbb{R}$ be a loss function and $\delta_n >0$. A function $u_{\theta^*} \in \mathcal{M}_n \subset U$ is said to be a \textit{quasi-minimizer} of $\mathcal{L}$ if
    \begin{equation}\label{Eq: Quasi-minimizer}
        \mathcal{L}(u_{\theta^*}) \leq \inf_{u_{\theta} \in \mathcal{M}_n} \mathcal{L}(u_{\theta}) + \delta_n.
    \end{equation}
\end{dfn}

\subsection{Error estimate for PINNs}

\begin{thm}\label{Thm: error estimate PINNs}
    Let $\delta_n > 0$ and $u_{\theta^*} \in \mathcal{M}_n$ be a quasi-minimizer of \cref{Eq: loss function PINNs} satisfying \cref{Eq: Quasi-minimizer}. If $u$ is the solution of \cref{Eq: steady state equation}, we have
    \begin{equation*}
        \|u-u_{\theta^*}\|_{U} \leq \left( \dfrac{1+\alpha}{1-\alpha} \right)\inf_{u_{\theta} \in \mathcal{M}_n}\|u-u_{\theta}\|_{U} + \dfrac{\delta_n}{1-\alpha}.
    \end{equation*}
\end{thm}
\begin{pf}
    See details of the proof in \ref{section: error estimates proof}.    
\end{pf}

\subsection{Error estimate for RVPINNs}
To obtain a reliable error bound for the RVPINNs method, we need the following assumption. 

\textbf{Local Fortin's condition}\footnote{Item 1 of this condition has been modified from Assumption 1 in \cite{rojas2024robust}, as we will derive an a priori error estimate for RVPINNs in terms of the infimum over  $\mathcal{M}_n^{\theta^*,R}$, which is a more natural form compared to the infimum over $\mathrm{span}(\mathcal{M}_n^{\theta^*,R})$ used in \cite{rojas2024robust}.}: Let $\delta_n > 0$ and $u_{\theta^*} \in \mathcal{M}_n$ be a quasi-minimizer of $\mathcal{L}_{\mathrm{RVPINNs}}$ satisfying \cref{Eq: Quasi-minimizer}. There exists $R > 0$ such that for all $\theta \in B(\theta^*,R)$, an operator $\Pi_{\theta}:V\to V_M$ and $\theta$-independent constant $C_{\Pi} > 0$, satisfying:
\begin{enumerate}
    \item $b(u_{\theta^*}-u_{\theta},v-\Pi_{\theta}v) = 0, \quad\forall v \in V,$
    \item $\|\Pi_{\theta}v\|_V \leq C_{\Pi}\|v\|_V, \quad \forall v\in V,$
\end{enumerate}
where $B(\theta^*,R)$ is an open ball of center $\theta^*$ and radius $R$, with respect to a given norm of $\mathbb{R}^n$. We denote $\mathcal{M}_n^{\theta^*,R} := \{ u_{\theta}\in \mathcal{M}_n : \theta \in B(\theta^*,R)\}$.

\begin{thm}\label{Thm: Error estimates variational problems}
    Let $\delta_n > 0$ and $u_{\theta^*} \in \mathcal{M}_n$ be a quasi-minimizer of \cref{Eq: loss function VPINNs sup} satisfying \cref{Eq: Quasi-minimizer}. If the local Fortin's condition is satisfied, we have
    \begin{equation*}
        \|u-u_{\theta^*}\|_{U} \leq \left( 1+\frac{2(1+\alpha)C_{\Pi}}{1-\alpha} \right)\inf_{u_{\theta} \in \mathcal{M}_n^{\theta^*,R}}\|u-u_{\theta}\|_{U} + \left( \frac{C_{\Pi}}{1-\alpha} \right)\delta_n.
    \end{equation*}
\end{thm}

\begin{pf}
    See details of the proof in \ref{section: error estimates proof}.    
\end{pf}

\section{Numerical examples}\label{sec: Numerical examples}

In this section, we show the performance of PINNs and RVPINNs for Neumann series problems by applying them to various 1D and 2D dynamical systems. We focus on systems where Perron-Frobenius operators are non-expansive under the $L^2$-norm. 

\subsection{1D dynamical systems}

We consider the dynamical system described by the tent map, which is a simple, continuous, piecewise-linear function exhibiting chaotic behavior. The tent map $S:[0,1] \to [0,1]$ is defined as:
\begin{equation*}
    S(x) = \left\{ \begin{matrix}
        2x, & x \in [0,\frac{1}{2}), \\
        2-2x, & x \in [\frac{1}{2},1].
    \end{matrix} \right.
\end{equation*}
By applying equation \cref{eq: Perron-Frobenius closed interval}, the Perron-Frobenius operator associated with the tent map can be derived in the form of
\begin{align}
    \mathcal{P}f(x) = \dfrac{1}{2}f\left(\dfrac{x}{2}\right) + \dfrac{1}{2}f\left(1-\dfrac{x}{2}\right) \label{Eq: PF of tent map},
\end{align}
and we can verify that it is non-expansive under the $L^2$-norm. Indeed, 
\begin{align*}
    \|\mathcal{P}f\|_{L^2}^2 &= \int_{0}^{1} \left( \dfrac{1}{2}f\left(\dfrac{x}{2}\right) + \dfrac{1}{2}f\left(1-\dfrac{x}{2}\right)\right)^2 dx \\
    &\leq \dfrac{1}{2}\int_0^1 f^2\left( \dfrac{x}{2} \right) + f^2\left( 1 - \dfrac{x}{2} \right) dx = \int_0^1 f^2(x) dx = \| f \|_{L^2}^2.
\end{align*}

Because it is challenging to obtain analytical solutions of \cref{Eq: steady state equation}, we design our examples in reverse order: We specify a desired solution and then determine the corresponding initial density. This approach allows us to create well-defined test cases.
For the Perron-Frobenius operator $\mathcal{P}$ given in \cref{Eq: PF of tent map}, we implement PINNs and RVPINNs to approximate two types of solutions:
\begin{enumerate}
    \item A smooth solution $u^{[1]}(x) = \exp(x)$,
    \item A singular solution $u^{[2]}(x) = 1 + x^{-1/3}$.
\end{enumerate}
To construct the initial densities $f_0^{[i]}$ corresponding to these solutions, we substitute the desired solution $u^{[i]}$ into \cref{Eq: steady state equation} and derive the following initial densities: 
\begin{enumerate}
    \item $f_0^{[1]}(x) = \exp(x) - \dfrac{\alpha}{2}\exp\left(\dfrac{x}{2}\right) - \dfrac{\alpha}{2}\exp\left(1-\dfrac{x}{2}\right)$,
    \item $f_0^{[2]}(x) = (1-\alpha) + x^{-1/3} - \dfrac{\alpha}{2}\left(\dfrac{x}{2}\right)^{-1/3} - \dfrac{\alpha}{2}\left(1-\dfrac{x}{2}\right)^{-1/3}$.
\end{enumerate}
We note that $f_0^{[1]}$ and $f_0^{[2]}$ are nonnegative when $\alpha \leq 2/(1+e)$ and $\alpha \leq 2/(1+\sqrt[3]{2})$, respectively.

\begin{figure}[t]
    \subfloat[PINNs approximation and solution]{\includegraphics[width=6.5cm]{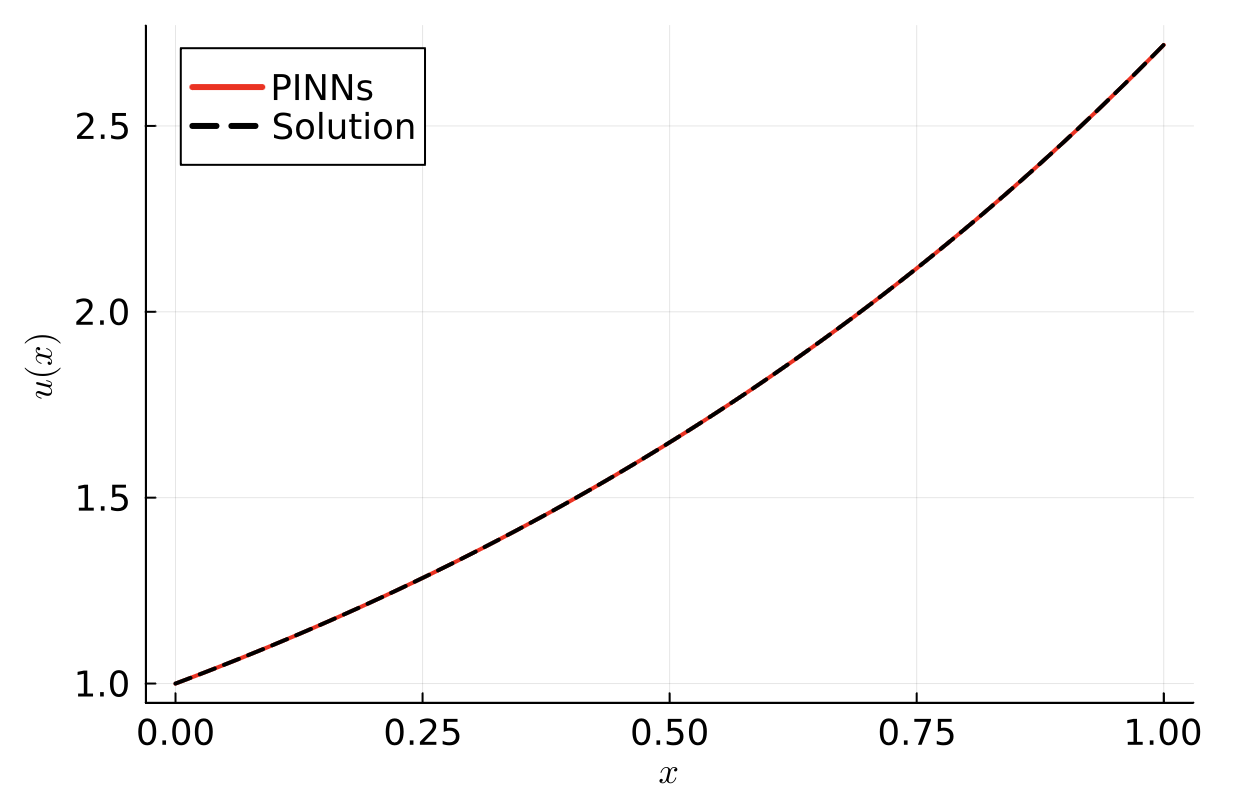}}\hfill 
    \subfloat[$L^2$-norm error for different number of neurons]{\includegraphics[width=6.5cm]{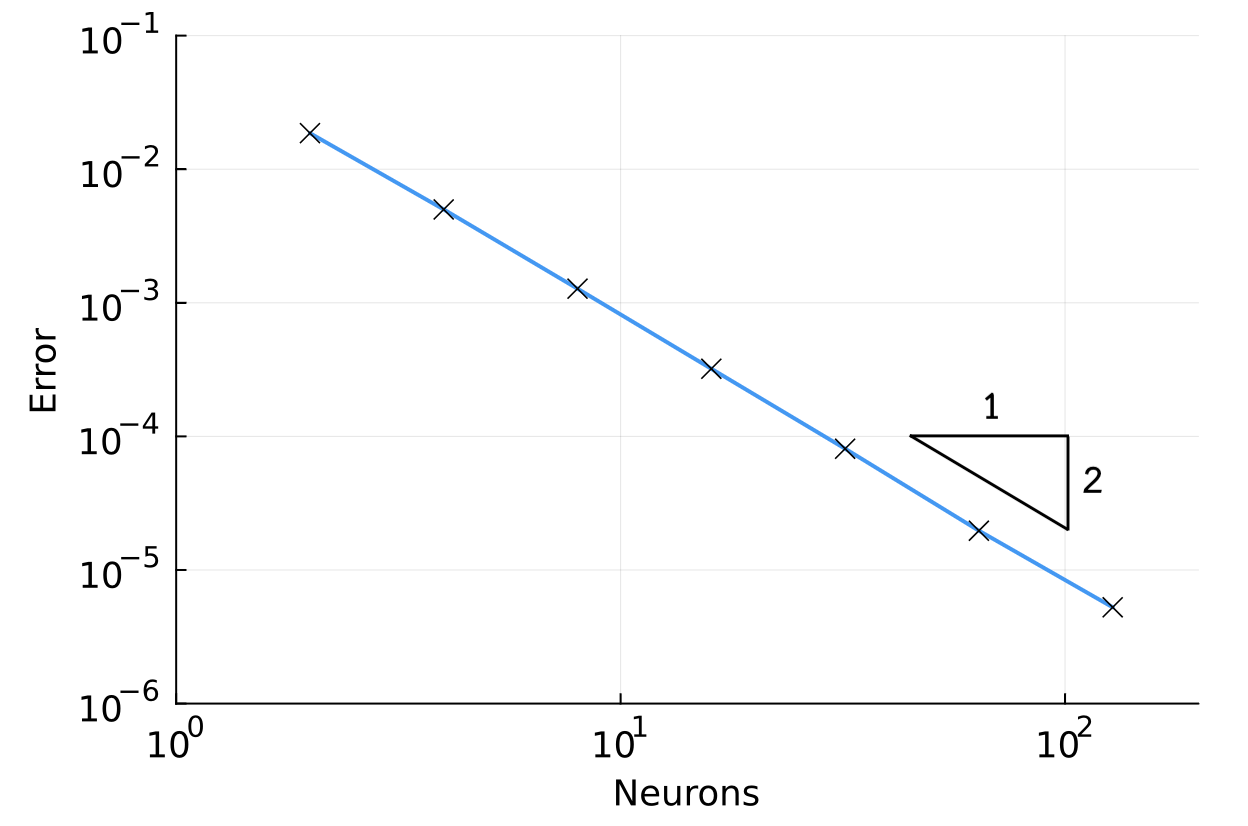}}
    \caption{PINNs approximation for $u^{[1]}$ with $\alpha = 0.5$.}
    \label{Figure: ANN approximation for tent map PINNs}
\end{figure}

To show the performance of the PINNs approach, we apply it to approximate the solution $u^{[1]}$. In the loss function \cref{Eq: loss function PINNs}, we use the corresponding initial density $f_0^{[1]}$, set the damping parameter to $\alpha = 0.5$, and use the Perron-Frobenius operator $\mathcal{P}$ as defined in \cref{Eq: PF of tent map}. The neural network architecture consists of a single hidden layer with $n$ neurons and the ReLU activation function ($\mathrm{ReLU}(x) = \max(x,0)$), i.e.,
\begin{equation}
     u_{\theta}(x) := \sum_{j=1}^{n} c_j\textrm{ReLU}(w_jx+b_j)+c_0,\label{Eq: M_n ReLU}
\end{equation}
where $c_j,w_j,b_j \in \mathbb{R}$ are trainable parameters.
The neural network is initialized to produce a piecewise-linear approximation over a uniform partition of $[0,1]$ into $n$ subintervals, with the outer parameters $c_j$ chosen to minimize the PINNs loss function. To approximate the integral in the loss function, we use a $101$-point Gauss-Kronrod quadrature rule. We train the neural network using the BFGS algorithm to minimize the loss function. 
Figure \ref{Figure: ANN approximation for tent map PINNs} (a) shows the PINNs approximation using the neural network with $n = 32$ neurons, which presents a good approximation to the solution. Figure \ref{Figure: ANN approximation for tent map PINNs} (b) demonstrates $L^2$-norm errors for trained neural networks, $\|u-u_\theta\|_{L^2}$, for different numbers of neurons.  The slope of the plot is approximately $-2$, confirming that the convergence rate is $O(n^{-2})$, as expected for continuous piecewise-linear approximations (see Remark \ref{remark: error bound piecewise linear approx}).

To compare PINNs with the fixed-grid-based method, we consider the singular solution $u^{[2]}$. In the fixed-grid-based method, the domain is partitioned into $n$ equal-length subintervals, and hat functions defined over these intervals are used as basis functions. In the PINNs approach, we use an $n$-neuron neural network as defined in \cref{Eq: M_n ReLU}. The inner parameters $w_j$ and $b_j$ are initialized to generate breakpoints at $r,r^2,\ldots,r^{n-1}$, and $0$ with $r = 0.662$, while the outer parameters $c_j$ are chosen to minimize the PINNs loss function for the given inner parameters. Under these settings, the fixed-grid-based method and PINNs produce piecewise-linear approximations with $n$ breakpoints. Since the target function is singular, we use a finer quadrature rule with $501$ points to improve integration accuracy. 

\begin{figure}[t]
    \subfloat[PINNs approximation, fixed-grid approximation and solution.]{\includegraphics[width=6.5cm]{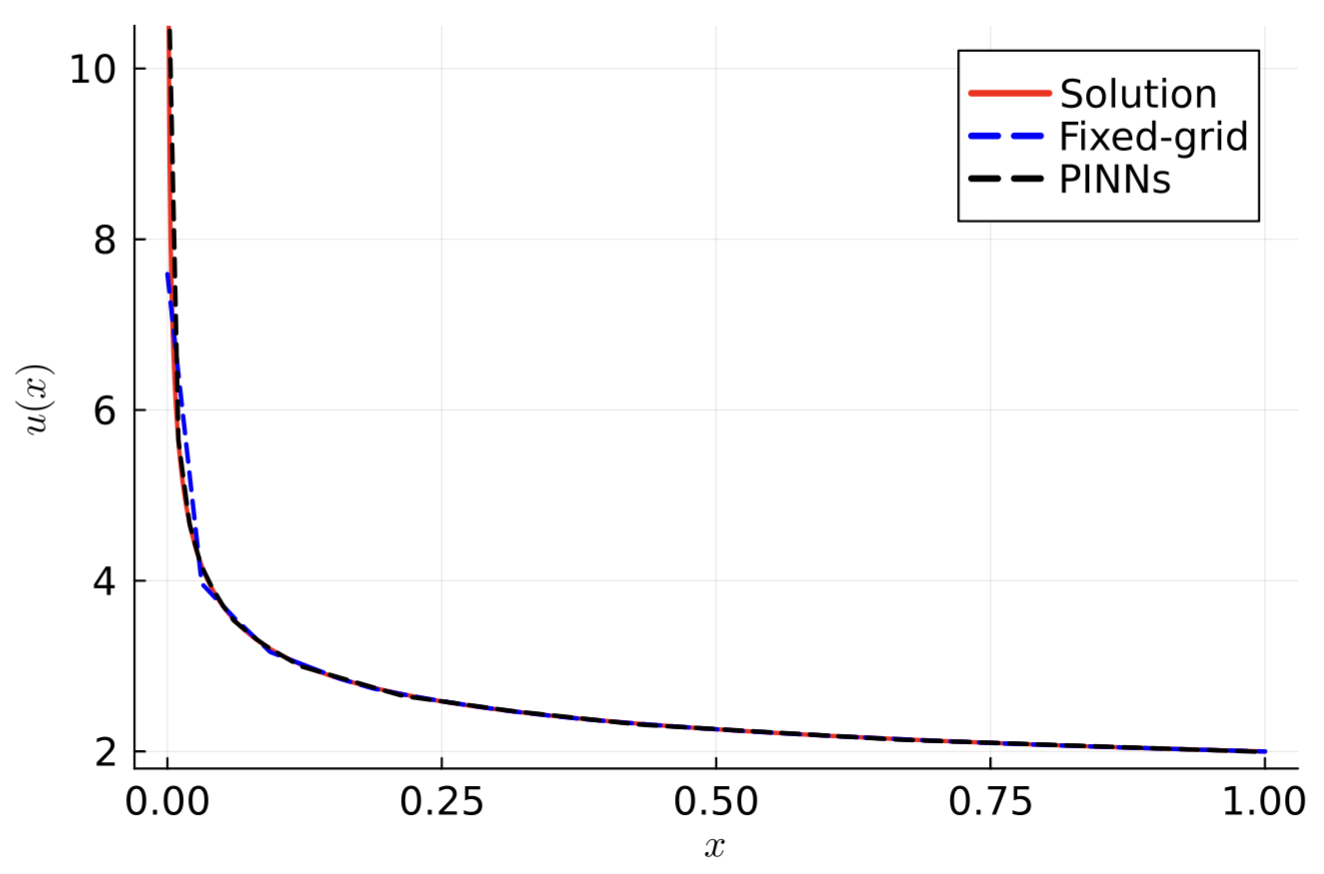}}\hfill 
    \subfloat[$L^2$-norm error comparison. The corresponding experimental orders of convergence are reported in Table \ref{Table: EOC fixed-grid PINNs}.]{\includegraphics[width=6.5cm]{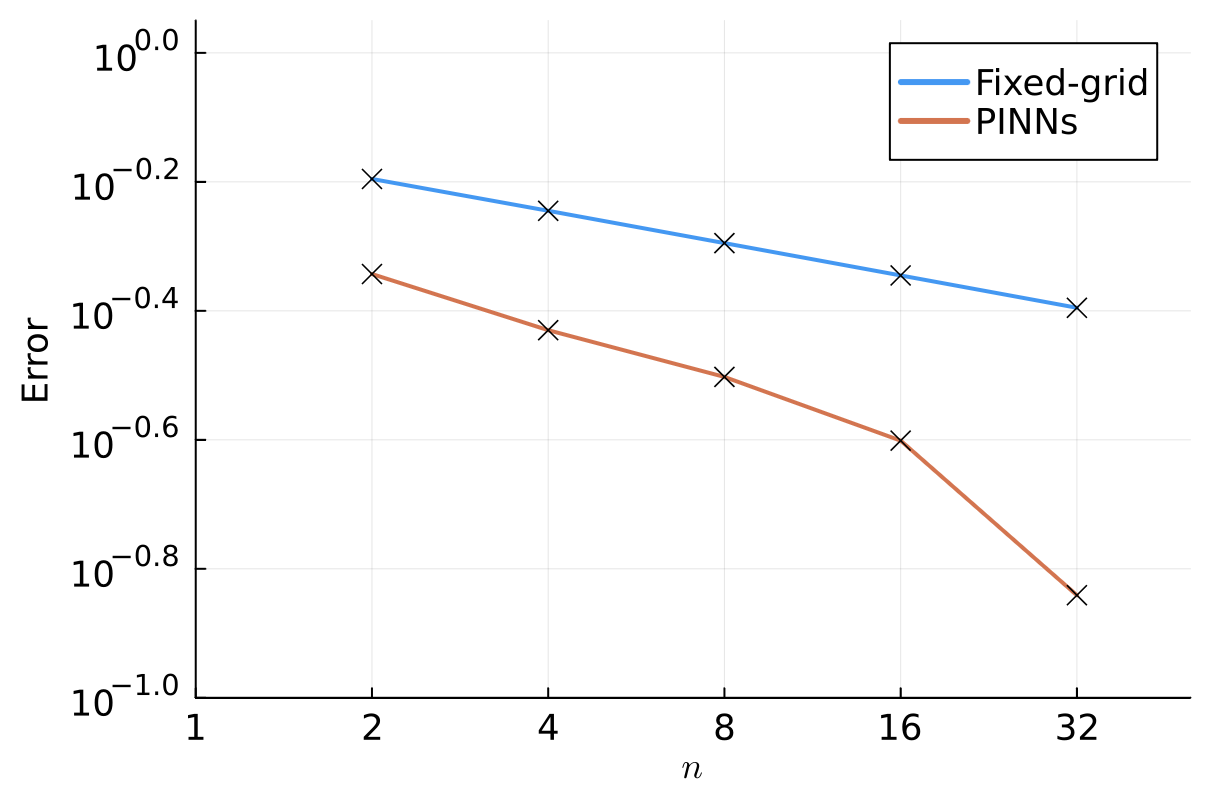}}
    \caption{PINNs approximation for $u^{[2]}$ with $\alpha = 0.5$.}
    \label{Figure: ANN approximation for tent map PINNs singular}
\end{figure}

\begin{table}[]
    \centering
    \begin{tabular}{  c  c  c }
        \hline
        $n$  & Fixed-grid & ~ PINNs ~ \\
        \hline
        $2 - 4$ & $0.164$ & $0.289$ \\
        $4 - 8$ & $0.166$ & $0.241$ \\
        $8 - 16$ & $0.167$ & $0.328$ \\
        $16 - 32$ & $0.167$ & $0.796$ \\
        \hline
    \end{tabular}
    \caption{Experimental orders of convergence for the fixed-grid-based method with $n$ breakpoints and PINNs with $n$ neurons.}
    \label{Table: EOC fixed-grid PINNs}
\end{table}

Figure \ref{Figure: ANN approximation for tent map PINNs singular} (a) shows that PINNs achieve greater accuracy than the fixed-grid-based method when both use $n=32$. Furthermore, Figure \ref{Figure: ANN approximation for tent map PINNs singular} (b) illustrates that the $L^2$-norm error for PINNs is smaller than that of the fixed-grid-based method, highlighting their superior approximation capability. Table \ref{Table: EOC fixed-grid PINNs} reports the experimental order of convergence (EOC) corresponding to Figure \ref{Figure: ANN approximation for tent map PINNs singular} (b). The EOC for the fixed-grid-based method is approximately $1/6$, consistent with the theoretical convergence rate $O(n^{-1/6})$ for a uniform grid, based on a Sobolev embedding argument as discussed in \cite[Example 5.8]{muga2019discrete}. Additionally, as expected for a neural network method, the EOC for PINNs shows a superior convergence rate.

\begin{figure}[t]
    \subfloat[RVPINNs approximation for $u^{[1]}$]{\includegraphics[width=6.5cm]{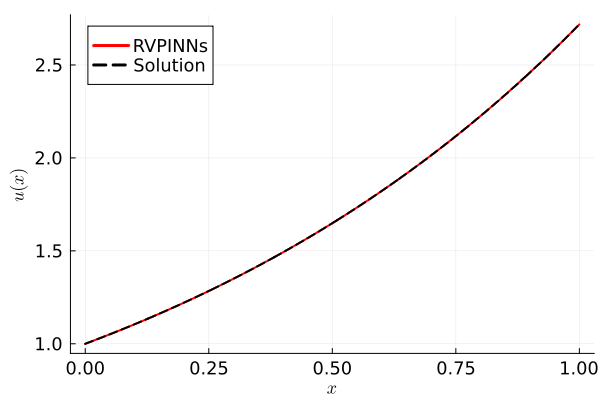}}\hfill 
    \subfloat[RVPINNs approximation for $u^{[2]}$]{\includegraphics[width=6.5cm]{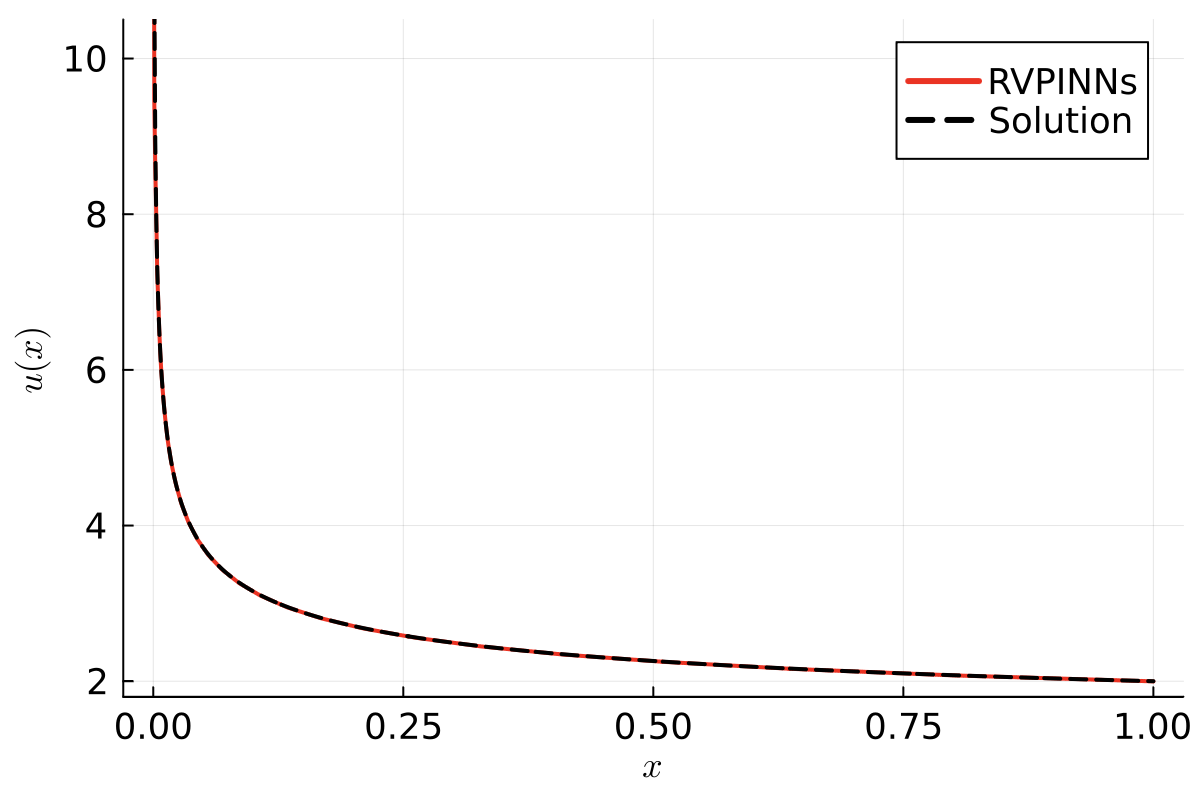}}
    \caption{RVPINNs approximation for the tent map example.}
    \label{Figure: ANN approximation for tent map VPINNs}
\end{figure}

The RVPINNs method provides an alternative approach for approximating solutions to Neumann series problems. We use the same neural network structure as the PINNs examples. The domain is partitioned into eight equal-length intervals, with the test functions chosen as normalized characteristic functions over these intervals. Figures \ref{Figure: ANN approximation for tent map VPINNs} (a) and \ref{Figure: ANN approximation for tent map VPINNs} 
(b) illustrate that the trained $32$-neuron neural networks, using the RVPINNs method, generate accurate approximations for both examples.

\begin{figure}[t]
    \centering
    \includegraphics[width=0.5\linewidth]{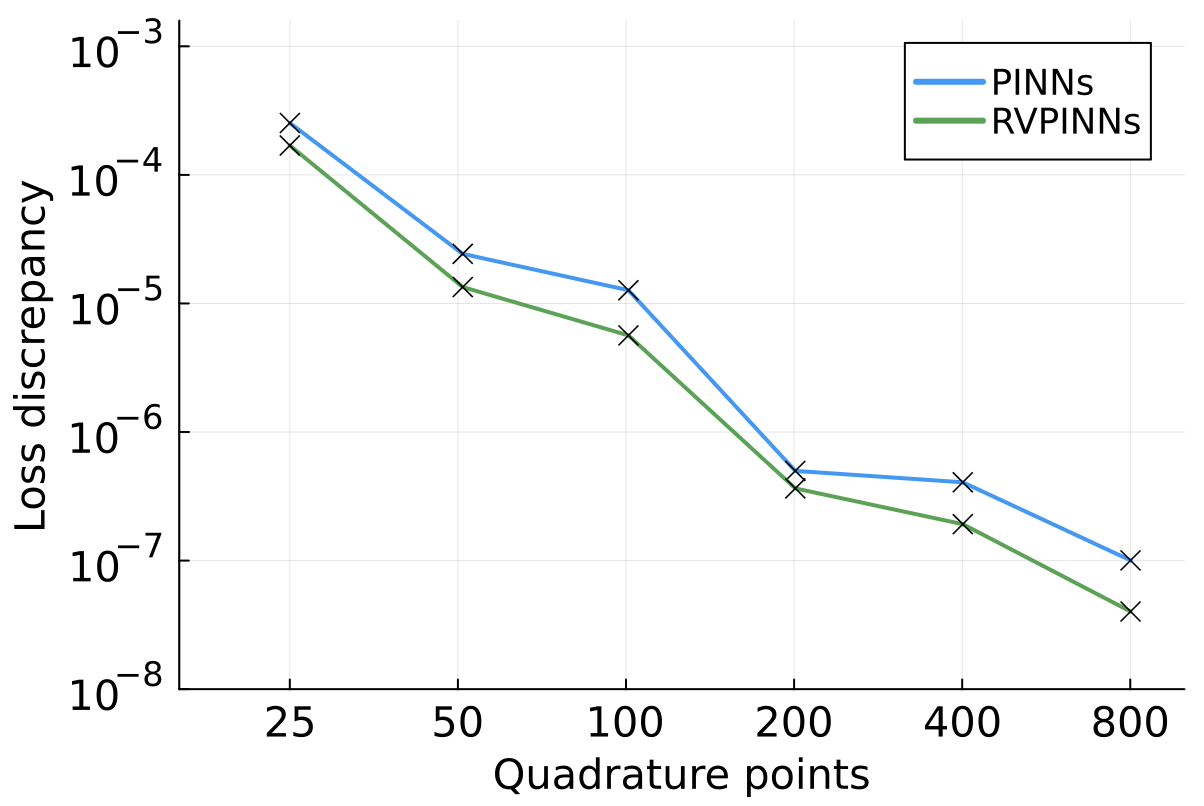}
    \caption{The discrepancy between the loss function computed using different numbers of quadrature points and the true value of the loss function.}
    \label{fig:difference_loss}
\end{figure}

Since our implementation relies on numerical integration using the Gauss-Kronrod quadrature rule, we examine the discrepancy between the loss function computed with $q$ quadrature points and that obtained using the adaptive Gauss-Kronrod quadrature rule, which we treat as the exact value (with an error tolerance below $10^{-8}$). As shown in Figure \ref{fig:difference_loss}, the approximation error in the loss function decreases as the number of quadrature points increases.

\subsection{2D dynamical systems}

In this section, we show the performance of PINNs and RVPINNs on 2D examples. We consider the boundary map in a circular domain and the standard map. We compare the numerical results with the $N$-term truncated sum, used as the ``exact'' Neumann series $u$ in \cref{Eq: power series P}. The $L^2$-norm of the truncation error is bounded by $\alpha^{N+1}\|f_0\|_{L^2}/(1-\alpha)$.

In a circular domain, the position and direction of rays can be described using coordinates $(\varphi, \psi)$, where 
$\varphi \in [0,2\pi)$ represents the position on the boundary, given by the polar angle, and $\psi \in (-\frac{\pi}{2},\frac{\pi}{2})$ denotes the direction of the reflected ray, measured as the angle between the inward-pointing normal vector and the outgoing ray (see Figure \ref{fig:circular domain}).
Define the phase space as $\Omega:=[0,2\pi) \times (-\frac{\pi}{2},\frac{\pi}{2})$. The boundary map $S_1:\Omega \to \Omega$ giving the coordinates of the next intersection is  defined by
\begin{equation*}
    S_1(\varphi,\psi)=(\varphi+\pi-2\psi \mod 2\pi,\psi),
\end{equation*}
and its inverse map is 
\begin{equation*}
    S_1^{-1}(\varphi,\psi)=(\varphi-\pi+2\psi \mod 2\pi,\psi).
\end{equation*}

\begin{figure}[t]
    \centering

\tikzset{every picture/.style={line width=0.75pt}} %set default line width to 0.75pt        

\begin{tikzpicture}[x=0.75pt,y=0.75pt,yscale=-1,xscale=1]

%Shape: Circle [id:dp959393517436786] 
\draw  [line width=1.5]  (100,95) .. controls (100,55.24) and (132.24,23) .. (172,23) .. controls (211.76,23) and (244,55.24) .. (244,95) .. controls (244,134.76) and (211.76,167) .. (172,167) .. controls (132.24,167) and (100,134.76) .. (100,95) -- cycle ;
%Straight Lines [id:da32067187439558986] 
\draw    (244,85) -- (190.35,26.47) ;
\draw [shift={(189,25)}, rotate = 47.49] [color={rgb, 255:red, 0; green, 0; blue, 0 }  ][line width=0.75]    (10.93,-3.29) .. controls (6.95,-1.4) and (3.31,-0.3) .. (0,0) .. controls (3.31,0.3) and (6.95,1.4) .. (10.93,3.29)   ;
\draw [shift={(244,85)}, rotate = 227.49] [color={rgb, 255:red, 0; green, 0; blue, 0 }  ][fill={rgb, 255:red, 0; green, 0; blue, 0 }  ][line width=0.75]      (0, 0) circle [x radius= 3.35, y radius= 3.35]   ;
%Straight Lines [id:da20417232067150826] 
\draw    (189,25) -- (124.93,42.47) ;
\draw [shift={(123,43)}, rotate = 344.74] [color={rgb, 255:red, 0; green, 0; blue, 0 }  ][line width=0.75]    (10.93,-3.29) .. controls (6.95,-1.4) and (3.31,-0.3) .. (0,0) .. controls (3.31,0.3) and (6.95,1.4) .. (10.93,3.29)   ;
\draw [shift={(189,25)}, rotate = 164.74] [color={rgb, 255:red, 0; green, 0; blue, 0 }  ][fill={rgb, 255:red, 0; green, 0; blue, 0 }  ][line width=0.75]      (0, 0) circle [x radius= 3.35, y radius= 3.35]   ;
%Straight Lines [id:da5540079650778997] 
\draw  [dash pattern={on 0.84pt off 2.51pt}]  (72,100) -- (270,102) ;
%Straight Lines [id:da36588971438616813] 
\draw  [dash pattern={on 4.5pt off 4.5pt}]  (171,101) -- (244,85) ;
%Straight Lines [id:da4965636478162765] 
\draw  [dash pattern={on 4.5pt off 4.5pt}]  (171,101) -- (189,25) ;
%Shape: Arc [id:dp6728766716449255] 
\draw  [draw opacity=0] (221.02,88.97) .. controls (221.01,88.65) and (221,88.32) .. (221,88) .. controls (221,80.01) and (224.13,72.75) .. (229.22,67.37) -- (251,88) -- cycle ; \draw  [color={rgb, 255:red, 144; green, 19; blue, 254 }  ,draw opacity=1 ] (221.02,88.97) .. controls (221.01,88.65) and (221,88.32) .. (221,88) .. controls (221,80.01) and (224.13,72.75) .. (229.22,67.37) ;  
%Shape: Arc [id:dp2918594131626562] 
\draw  [draw opacity=0] (205.68,93.18) .. controls (206.48,95.76) and (206.91,98.4) .. (206.99,101.03) -- (177,102) -- cycle ; \draw  [color={rgb, 255:red, 144; green, 19; blue, 254 }  ,draw opacity=1 ] (205.68,93.18) .. controls (206.48,95.76) and (206.91,98.4) .. (206.99,101.03) ;  
%Shape: Arc [id:dp4192582475938966] 
\draw  [draw opacity=0] (176.63,78.99) .. controls (183.75,80.86) and (190.2,85.34) .. (194.42,92.07) .. controls (196.37,95.19) and (197.68,98.52) .. (198.39,101.92) -- (169,108) -- cycle ; \draw  [color={rgb, 255:red, 208; green, 2; blue, 27 }  ,draw opacity=1 ] (176.63,78.99) .. controls (183.75,80.86) and (190.2,85.34) .. (194.42,92.07) .. controls (196.37,95.19) and (197.68,98.52) .. (198.39,101.92) ;  
%Shape: Arc [id:dp045569316825823325] 
\draw  [draw opacity=0] (183.31,52.64) .. controls (174.89,49.07) and (168.43,41.76) .. (166.02,32.8) -- (195,25) -- cycle ; \draw  [color={rgb, 255:red, 208; green, 2; blue, 27 }  ,draw opacity=1 ] (183.31,52.64) .. controls (174.89,49.07) and (168.43,41.76) .. (166.02,32.8) ;  
%Straight Lines [id:da9193074506316568] 
\draw  [dash pattern={on 0.84pt off 2.51pt}]  (170,2) -- (170,187) ;

% Text Node
\draw (195,98) node [anchor=north west][inner sep=0.75pt]  [color={rgb, 255:red, 144; green, 19; blue, 254 }  ,opacity=1 ] [align=left] {$\displaystyle \varphi $};
% Text Node
\draw (122,148) node [anchor=north west][inner sep=0.75pt]   [align=left] {};
% Text Node
\draw (222,71) node [anchor=north west][inner sep=0.75pt]  [color={rgb, 255:red, 144; green, 19; blue, 254 }  ,opacity=1 ] [align=left] {$\displaystyle \textcolor[rgb]{0.56,0.07,1}{\psi }$};
% Text Node
\draw (170,28) node [anchor=north west][inner sep=0.75pt]  [color={rgb, 255:red, 208; green, 2; blue, 27 }  ,opacity=1 ] [align=left] {$\displaystyle \textcolor[rgb]{0.82,0.01,0.11}{\psi '}$};
% Text Node
\draw (175,78) node [anchor=north west][inner sep=0.75pt]  [color={rgb, 255:red, 208; green, 2; blue, 27 }  ,opacity=1 ] [align=left] {$\displaystyle \textcolor[rgb]{0.82,0.01,0.11}{\varphi '}$};
% Text Node
\draw (249,74) node [anchor=north west][inner sep=0.75pt]  [color={rgb, 255:red, 144; green, 19; blue, 254 }  ,opacity=1 ] [align=left] {$\displaystyle ( \varphi ,\psi )$};
% Text Node
\draw (198,7) node [anchor=north west][inner sep=0.75pt]  [color={rgb, 255:red, 208; green, 2; blue, 27 }  ,opacity=1 ] [align=left] {$\displaystyle ( \varphi ',\psi ')$};

\end{tikzpicture}

\caption{Diagram of a ray trajectory in a circular domain. The boundary map $S_1$ sends $(\varphi,\psi)$ to the new ray coordinates $(\varphi',\psi')$.}
\label{fig:circular domain}
\end{figure}

The standard map is a mathematical model used to study chaotic behavior in dynamical systems. It is a discrete-time dynamical system that exhibits a transition from regular to chaotic motion, making it a popular example in the study of chaos theory.
The standard map $S_2$ is a mapping from $[0,2\pi)\times[0,2\pi)$ to itself defined by 
    \begin{equation*}
        S_2(\theta,p) = (\theta + p + K\sin(\theta) \mod 2\pi, p + K\sin(\theta) \mod 2\pi),
    \end{equation*}
and its inverse map is
    \begin{equation*}
        S_2^{-1}(\theta,p) = (\theta - p \mod 2\pi,p - K\sin(\theta - p) \mod 2\pi),
    \end{equation*}
where $K$ is a parameter that controls the amount of chaos in the system. In this work, we use $K=2.4$.

Note that the boundary map and the standard map satisfy Remark \ref{Remark: S diffeomorphism} with $\left|J_{S_i^{-1}}\right| = 1$, so the Perron-Frobenius operators associated with $S_i$ are in the form $\mathcal{P}f = f(S_i^{-1})$. 

\begin{figure}[p]
    \subfloat[2D plot of the truncated sum with $1{,}000$ terms]{\includegraphics[width=6cm]{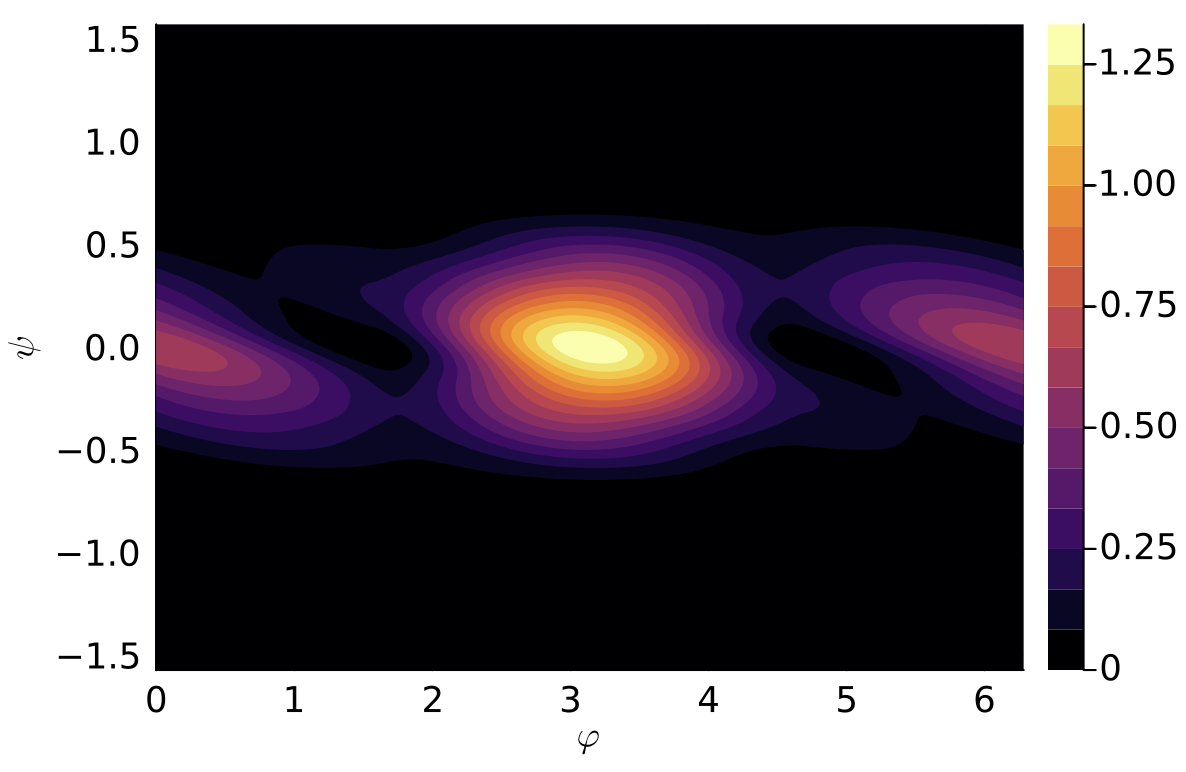}}\hfill 
    \subfloat[3D plot of the truncated sum with $1{,}000$ terms]{\includegraphics[width=6cm]{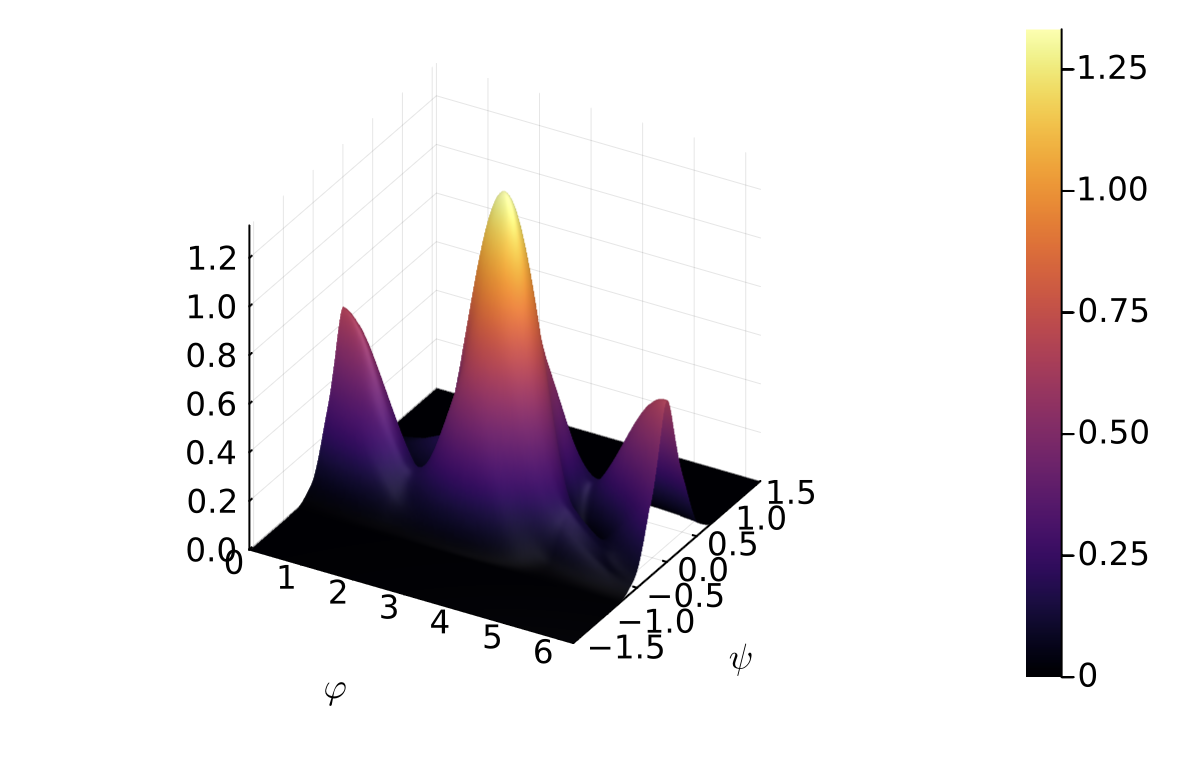}} \\
    \subfloat[2D plot of PINNs approximation]{\includegraphics[width=6cm]{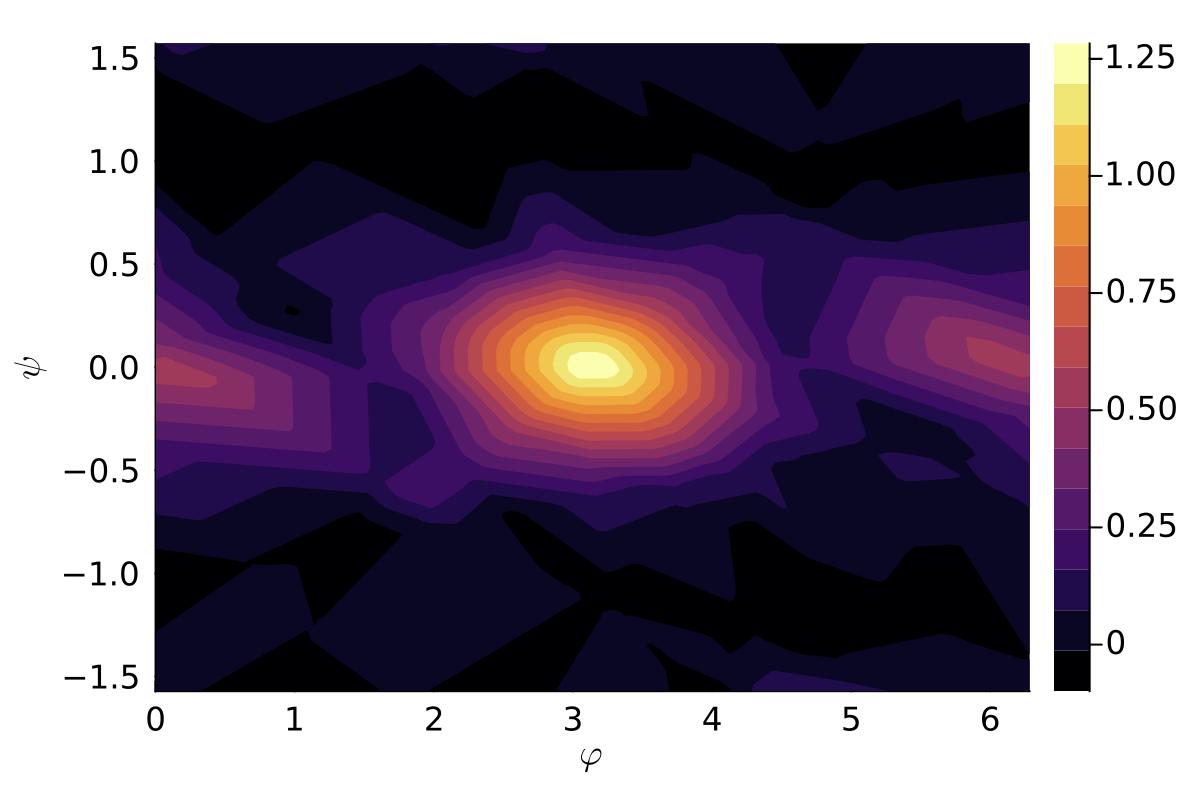}}\hfill 
    \subfloat[3D plot of PINNs approximation]{\includegraphics[width=6cm]{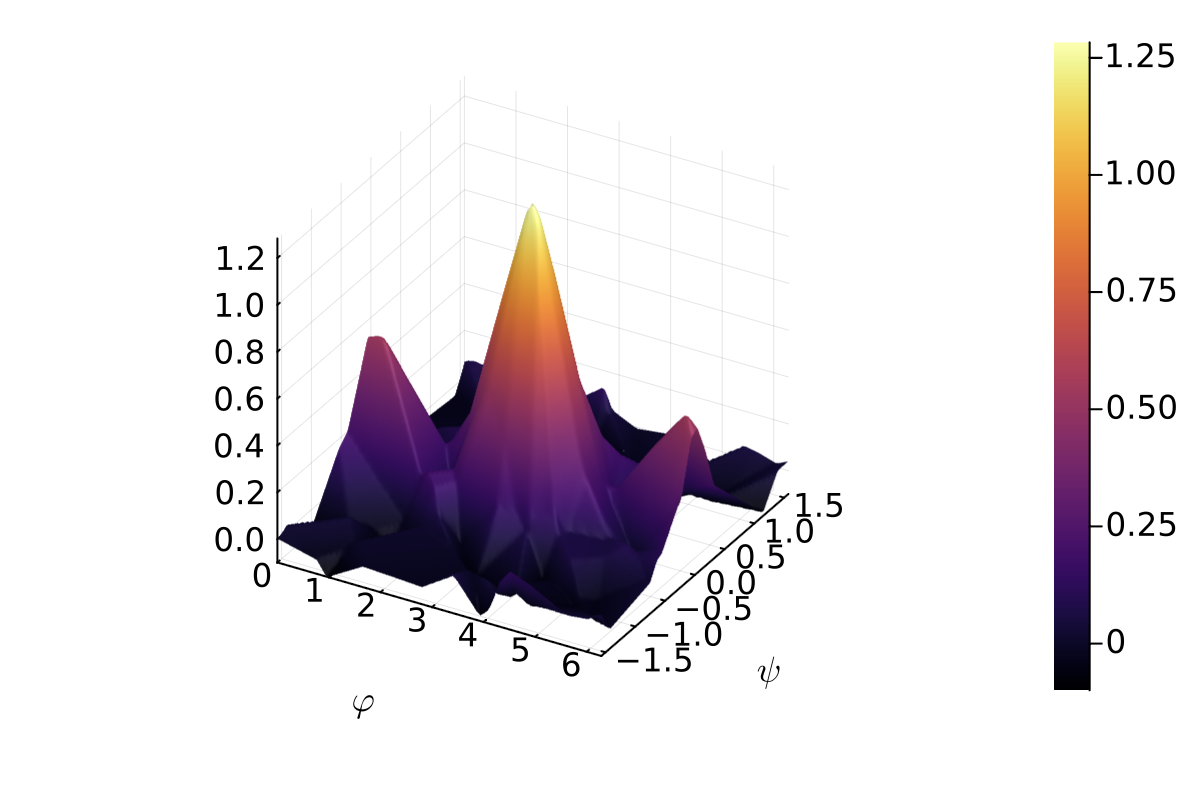}}\\
    \subfloat[2D plot of PINNs approximation with straight lines representing the affine transformations in the hidden layer]{\includegraphics[width=6cm]{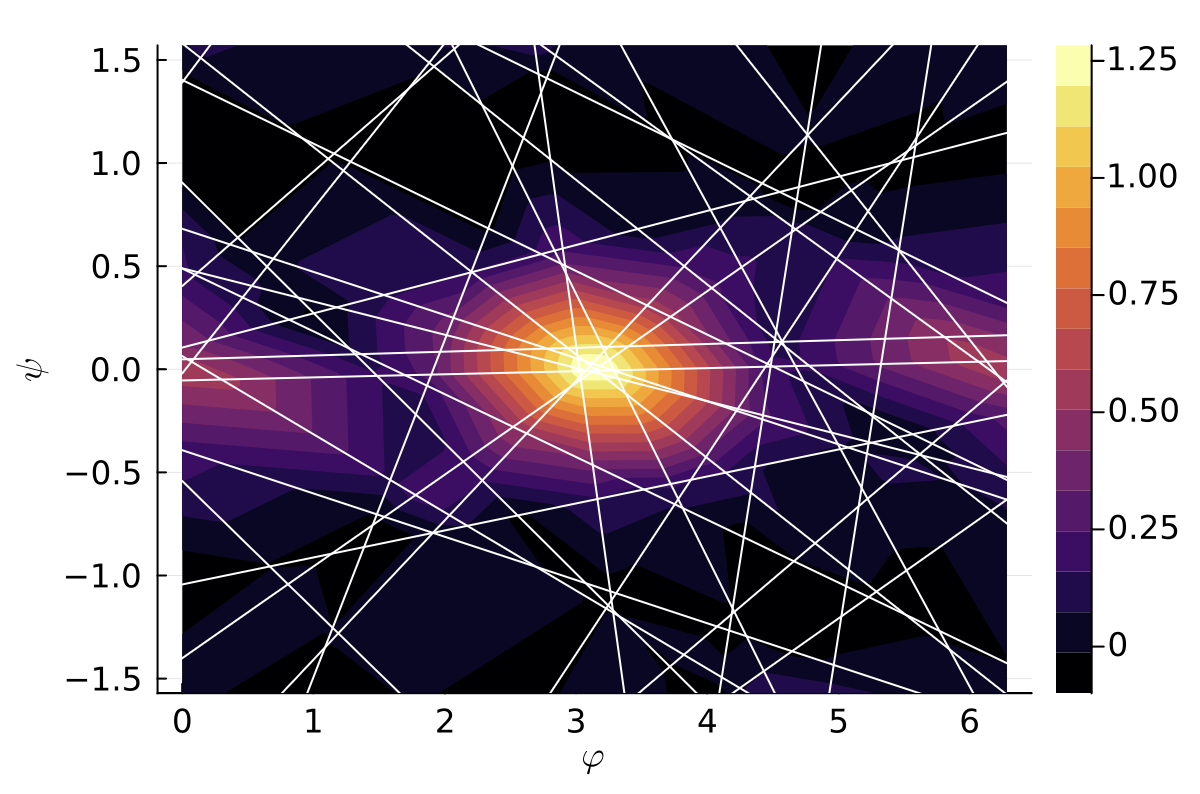}}\hfill 
    \subfloat[Loss plot]{\includegraphics[width=6cm]{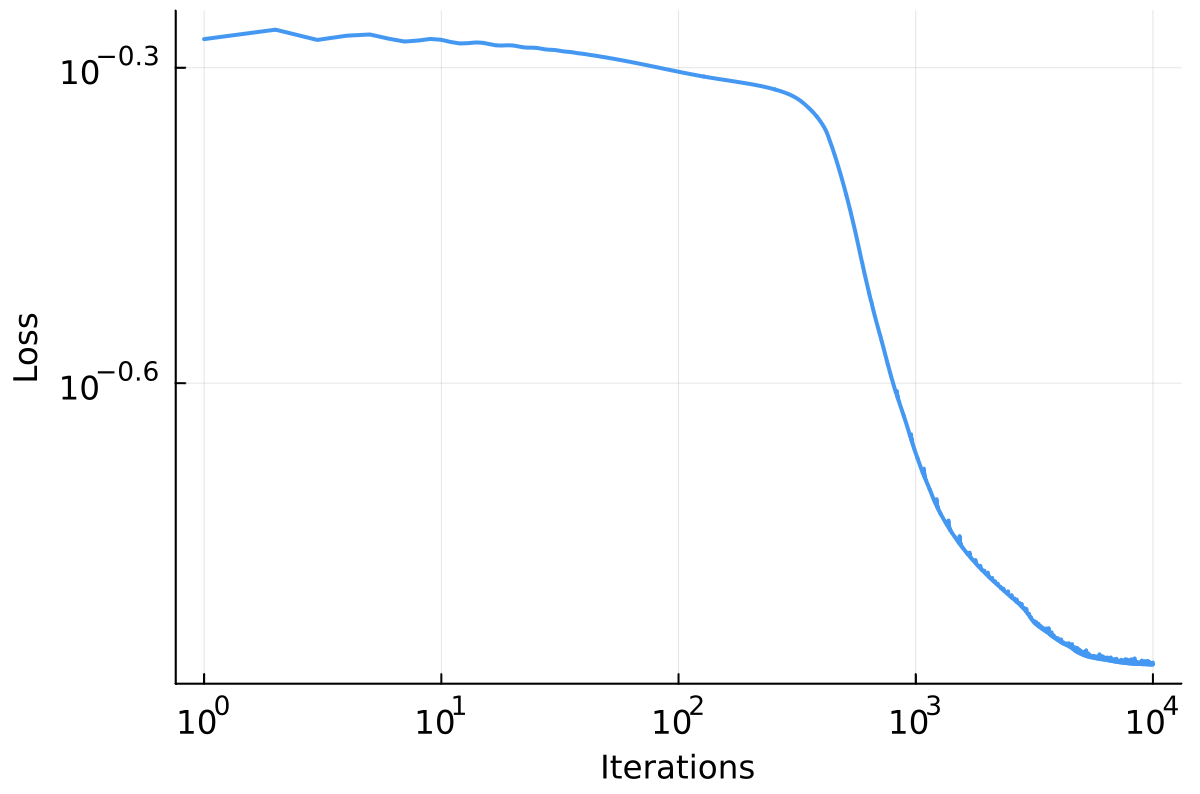}}
    \caption{PINNs approximation for the boundary map example.}
    \label{Figure: ANN approximation for phase map in circular domain}
\end{figure}

\begin{figure}[h] 
    \subfloat[2D plot of RVPINNs approximation]{\includegraphics[width=6cm]{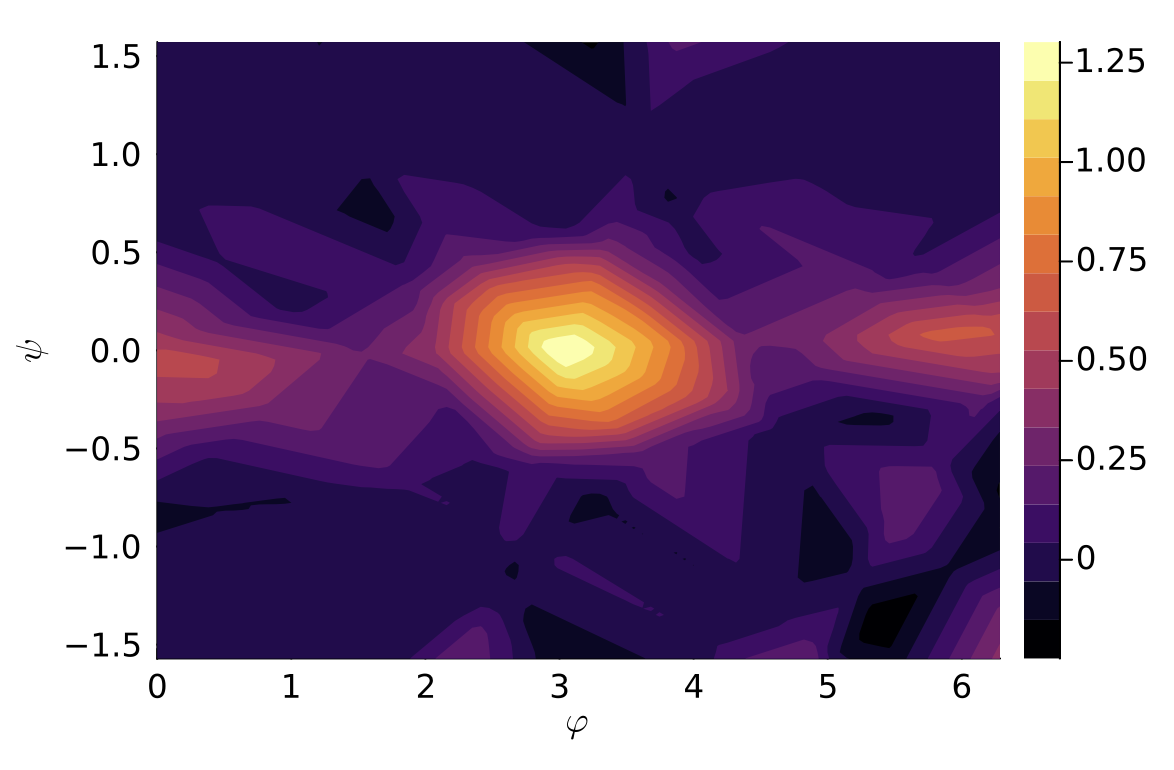}}\hfill 
    \subfloat[3D plot of RVPINNs approximation]{\includegraphics[width=6cm]{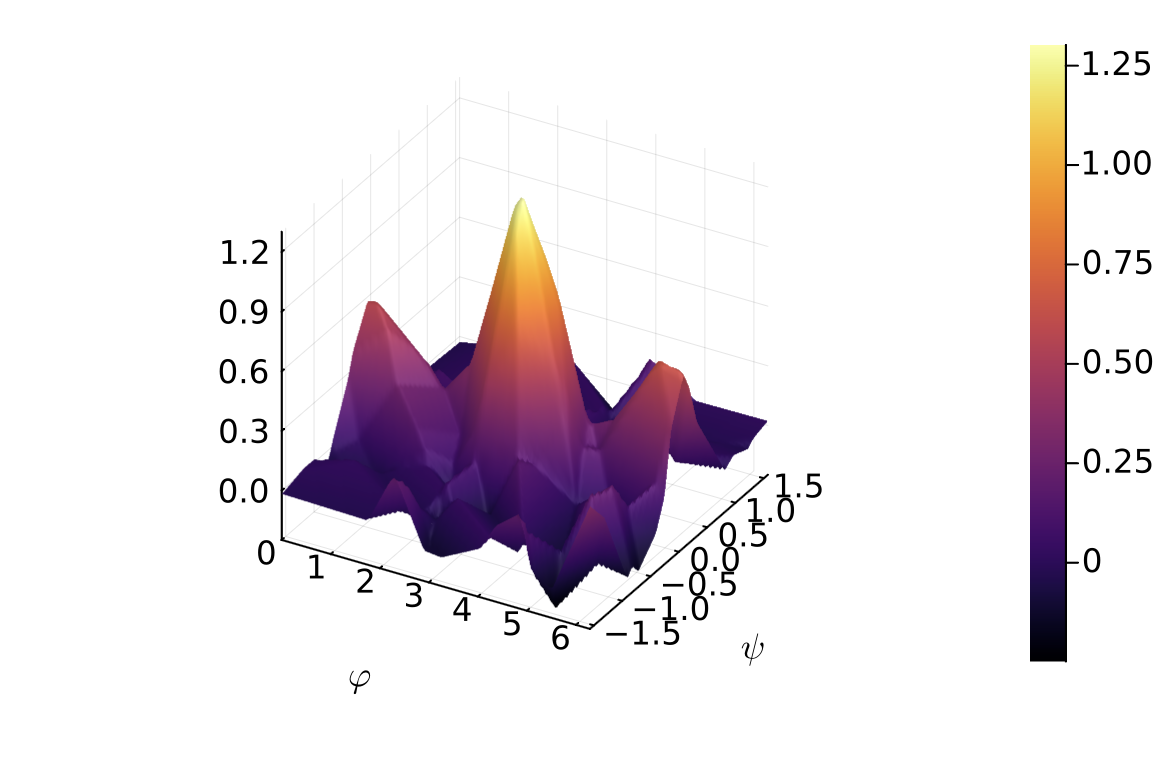}}\\
    \subfloat[2D plot of RVPINNs approximation with straight lines representing the affine transformations in the hidden layer]{\includegraphics[width=6cm]{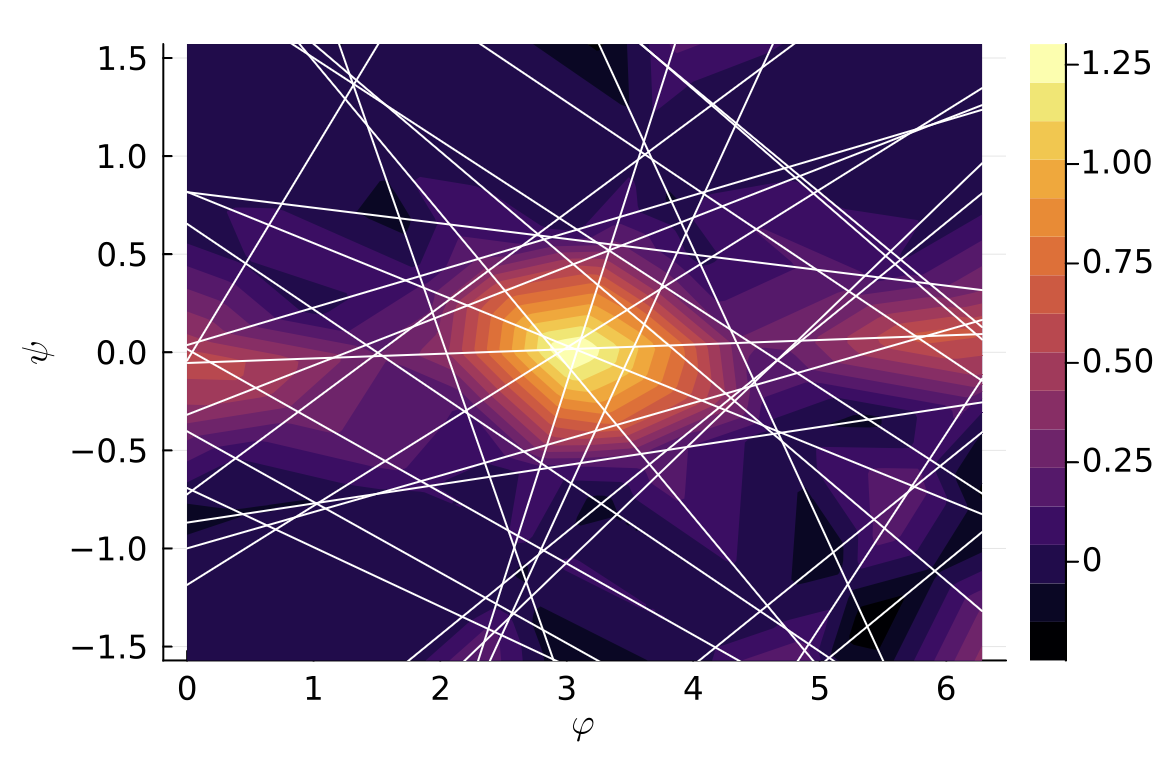}}\hfill 
    \subfloat[Loss plot]{\includegraphics[width=6cm]{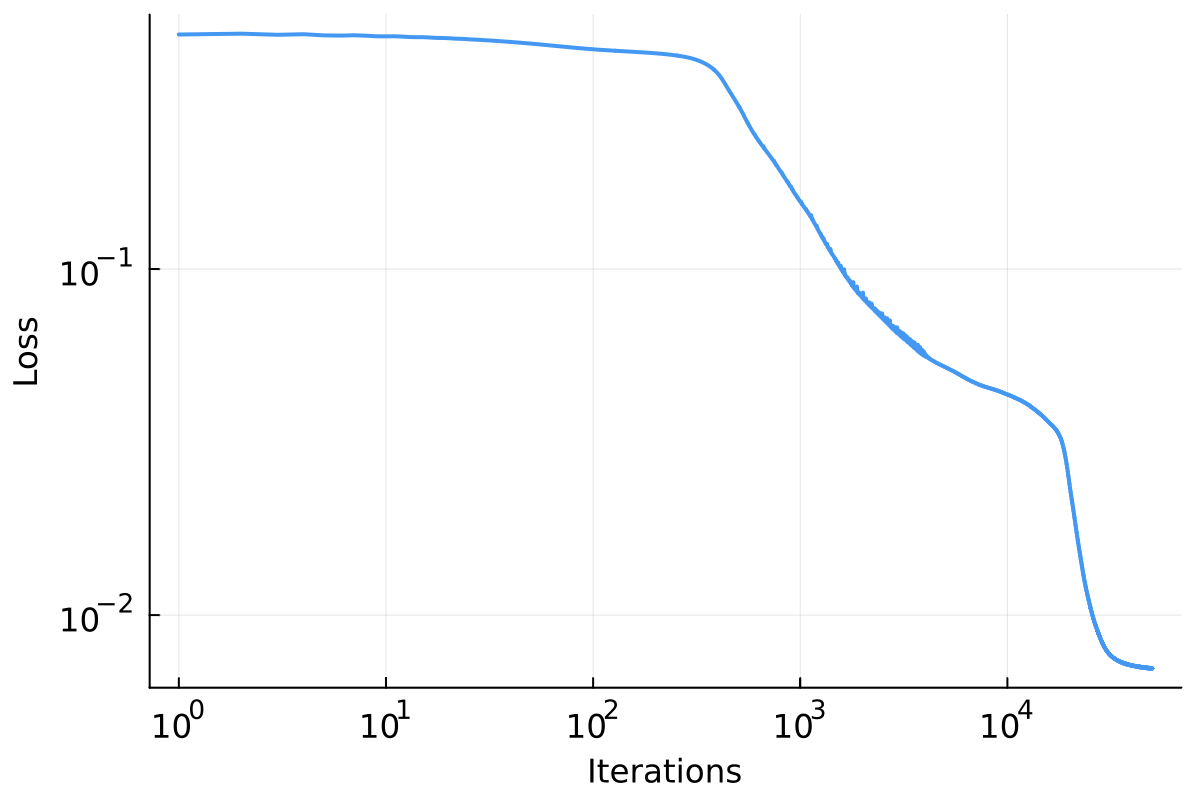}}
    \caption{RVPINNs approximation for the boundary map example.}
    \label{Figure: VPINNs steady state circular}
\end{figure}

We implement PINNs and RVPINNs to approximate the solution of \cref{Eq: steady state equation} and \cref{steady state variational problem}, respectively, for the Perron-Frobenius operator associated with the boundary map.
We set the damping parameter $\alpha = 0.5$. The initial density is given by $f_0(\varphi,\psi)=\cos^2(\varphi)\cos^2(2\psi)$ for $\pi/2 < \varphi < 3\pi/2$ and $-\pi/4 < \psi < \pi/4$, and $0$ elsewhere. Since $f_0$ is bounded and compactly supported, we have $f_0 \in L^2(\Omega)$. For RVPINNs, we partition the domain into $8 \times 8$ equal-sized cells, denoted as $\omega_1,\ldots,\omega_{64}$, and define $g_m$ as the normalized characteristic function on $\omega_m$.  We use a two-layer neural network, $32$ neurons, and the ReLU activation function, i.e.,
\begin{equation*}
     u_{\theta}(x,y) := \sum_{j=1}^{32} c_j\textrm{ReLU}(w_j^{(1)}x+w_j^{(2)}y+b_j)+c_0,
\end{equation*}
where $c_j,w_j^{(1)},w_j^{(2)},b_j \in \mathbb{R}$ are trainable parameters.
We approximate double integrals in PINNs and RVPINNs loss functions using a $101$-point Gauss-Kronrod quadrature rule for each variable. After training the neural network using Adam optimizer, we obtain approximations that closely match the truncated sum of \cref{Eq: power series P} with $1{,}000$ terms (``exact" solution), as shown in Figures \ref{Figure: ANN approximation for phase map in circular domain} and \ref{Figure: VPINNs steady state circular}. Despite using a relatively small neural network, it achieves a good approximation.

\begin{figure}[h]
    \subfloat[Classical phase space]{\includegraphics[width=6cm]{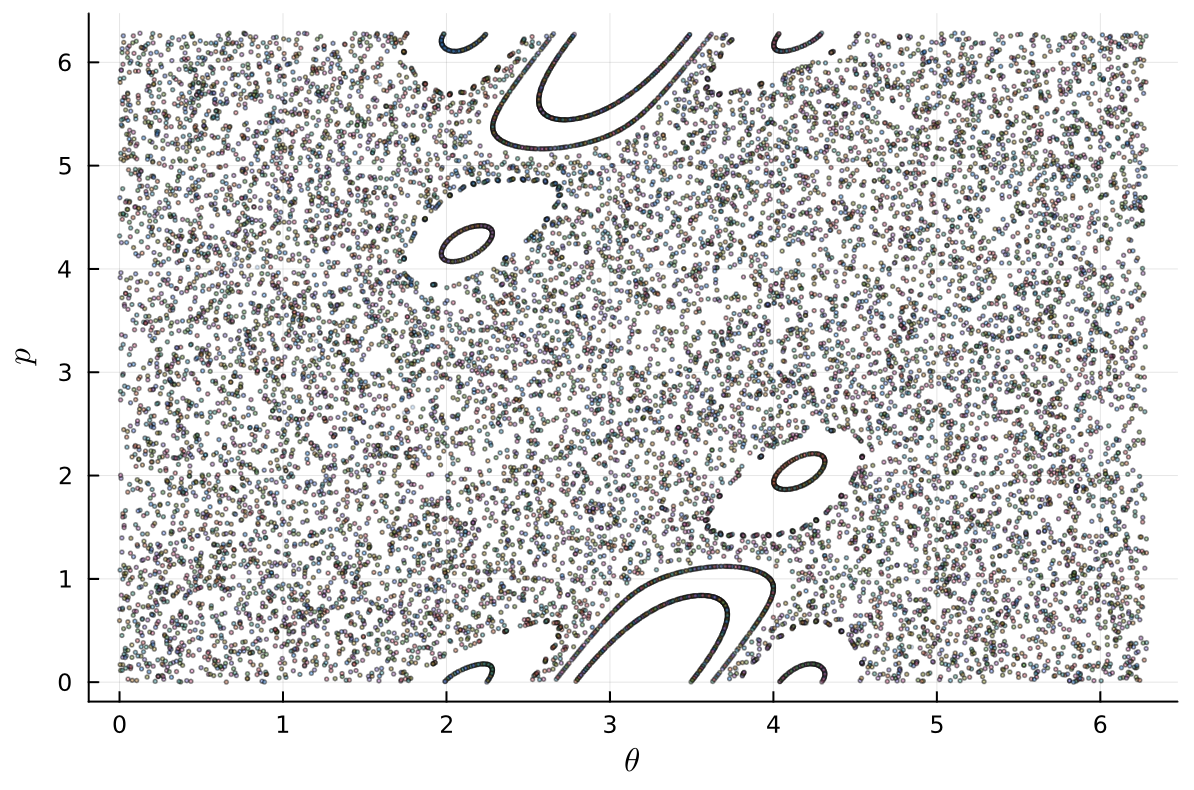}}\hfill 
    \subfloat[Initial density]{\includegraphics[width=6cm]{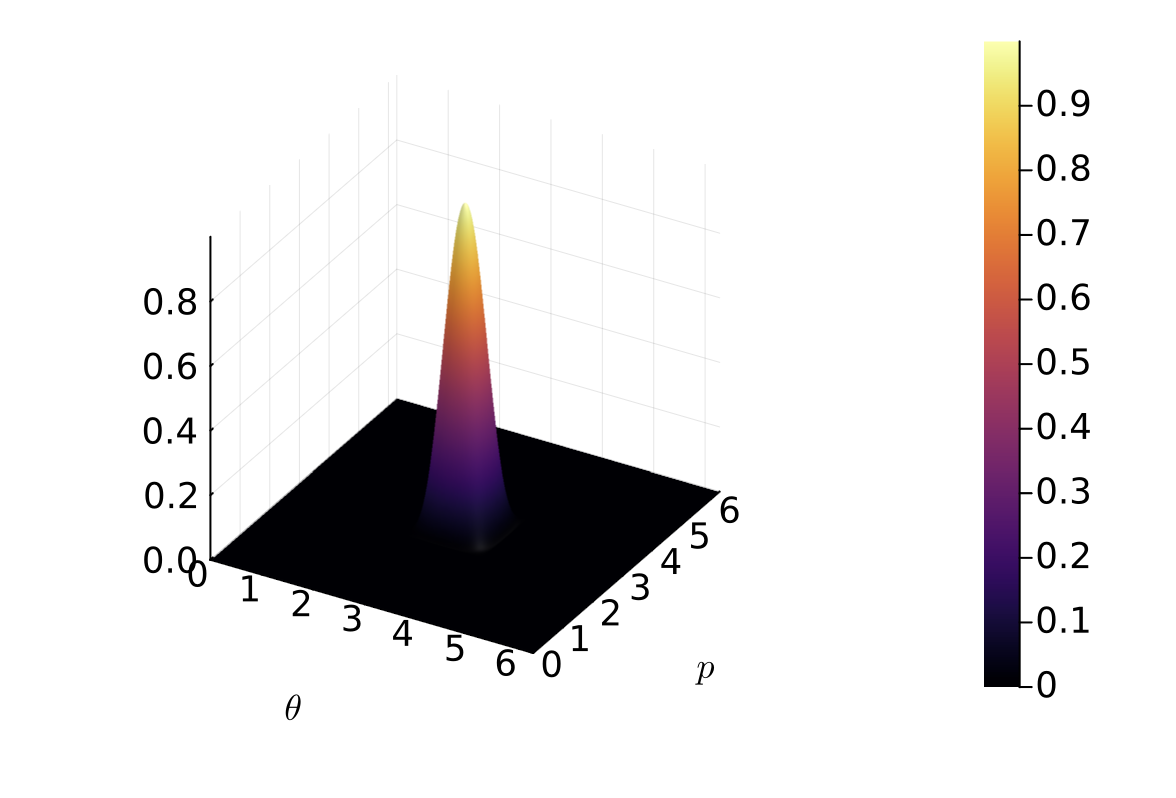}}\\
    \subfloat[$\alpha = 0.5$]{\includegraphics[width=6cm]{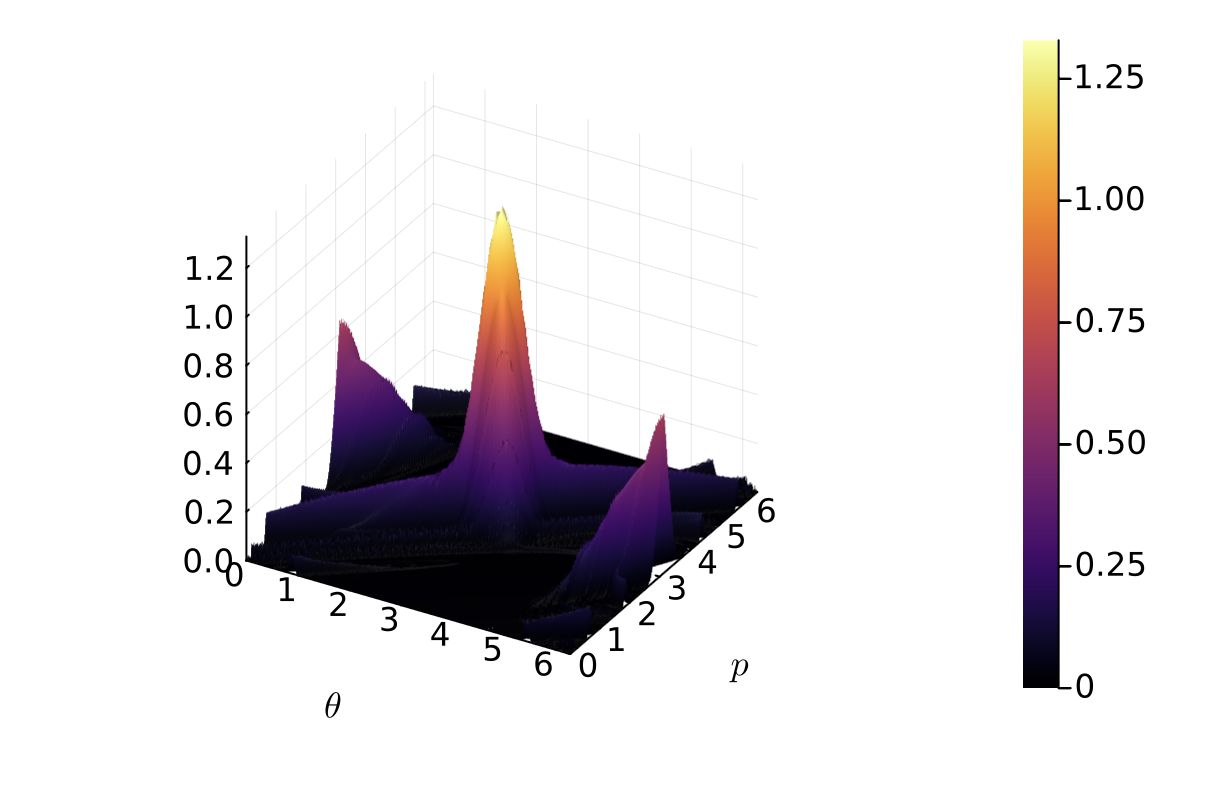}}\hfill 
    \subfloat[$\alpha = 0.9$]{\includegraphics[width=6cm]{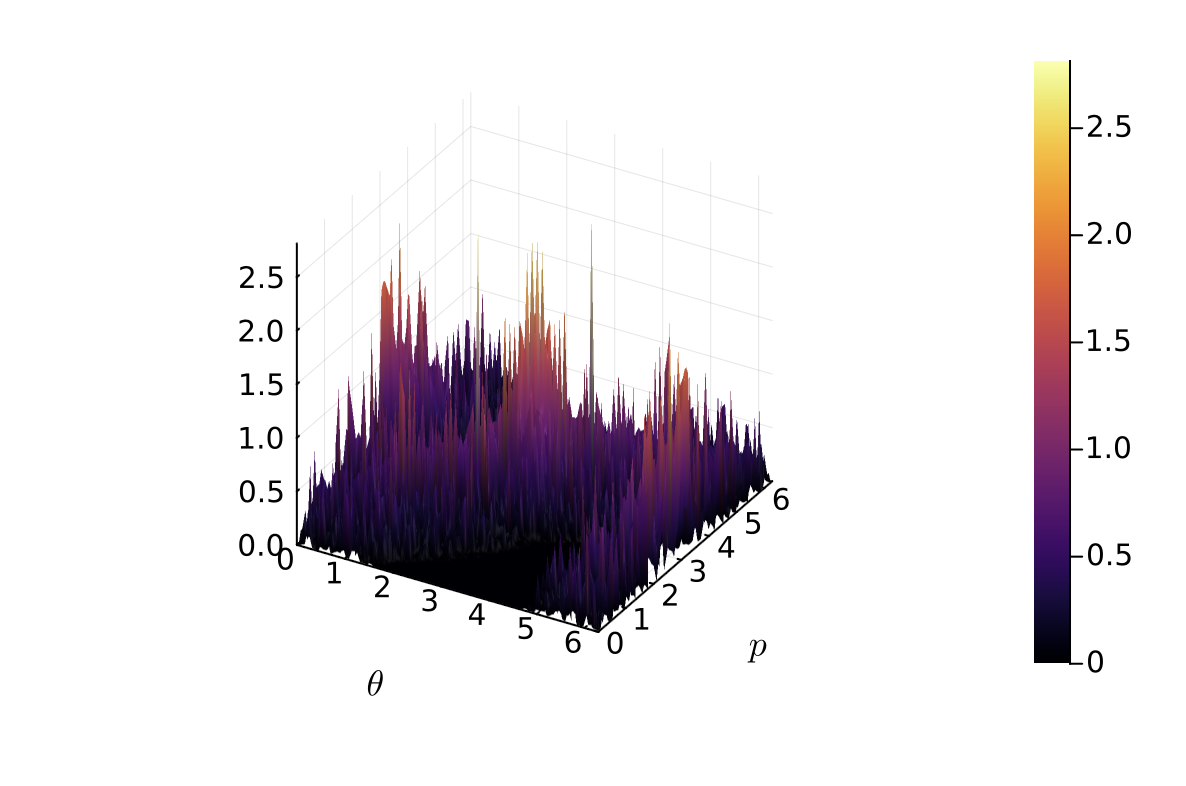}}
    \caption{(a) A classical phase space showing the distribution of trajectories of the standard map. (b) Initial density. (c) - (d) The truncated sums with $1{,}000$ terms using $\alpha = 0.5$ and $\alpha = 0.9$, respectively. }
    \label{Figure: VPINNs steady state standard map}
\end{figure}

\begin{figure}[p]
    \subfloat[Truncated sum]{\includegraphics[width=6cm]{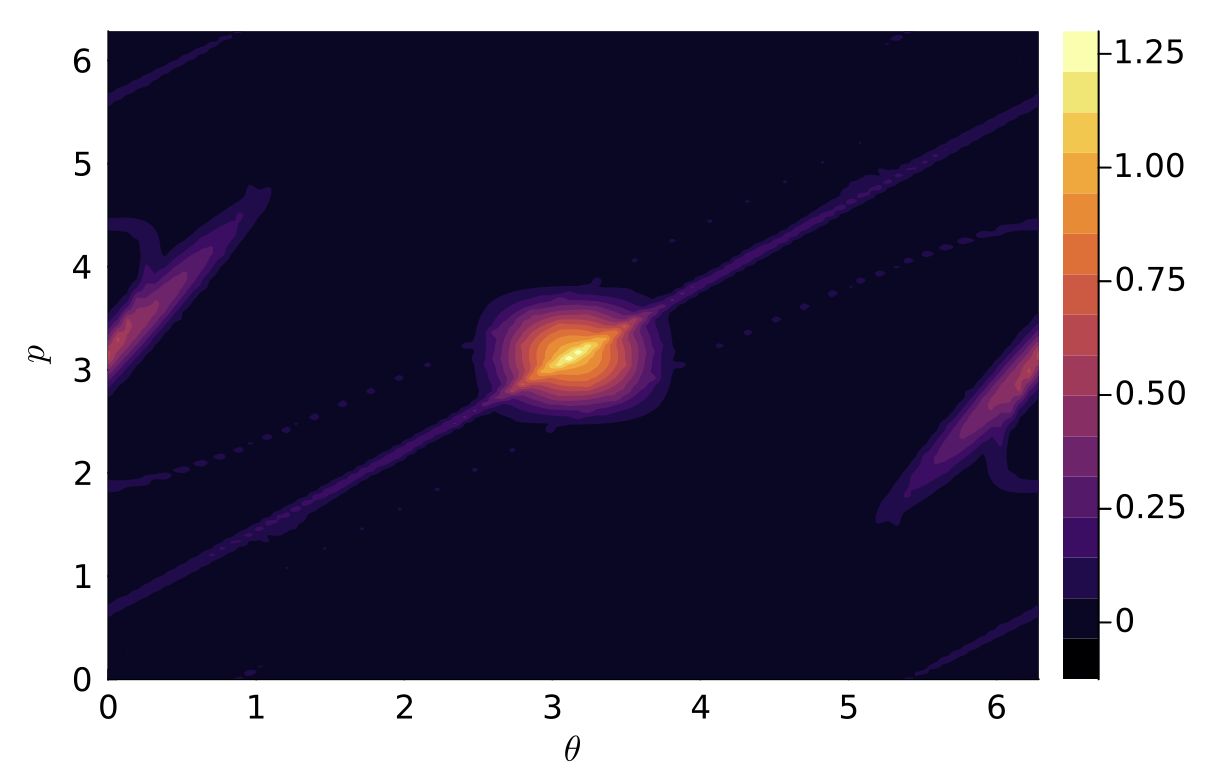}}\hfill
    \subfloat[Truncated sum]{\includegraphics[width=6cm]{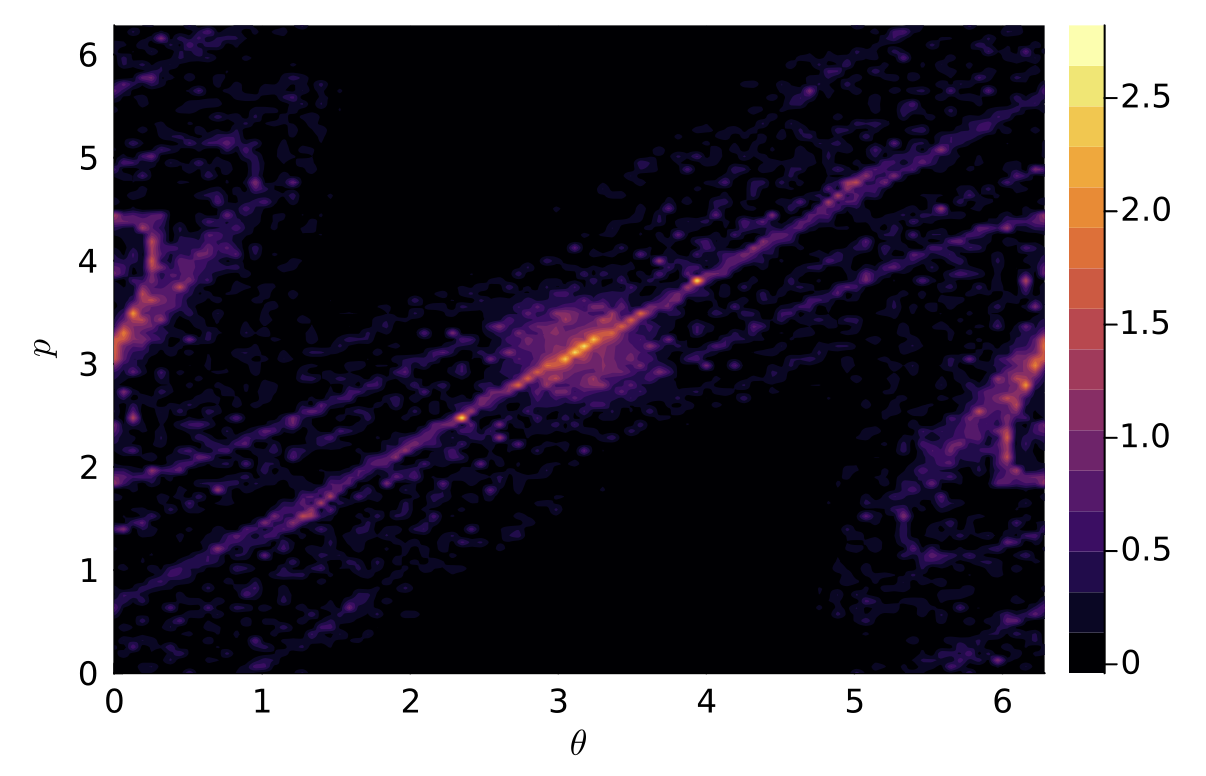}} \\
    \subfloat[Two-layer neural network]{\includegraphics[width=6cm]{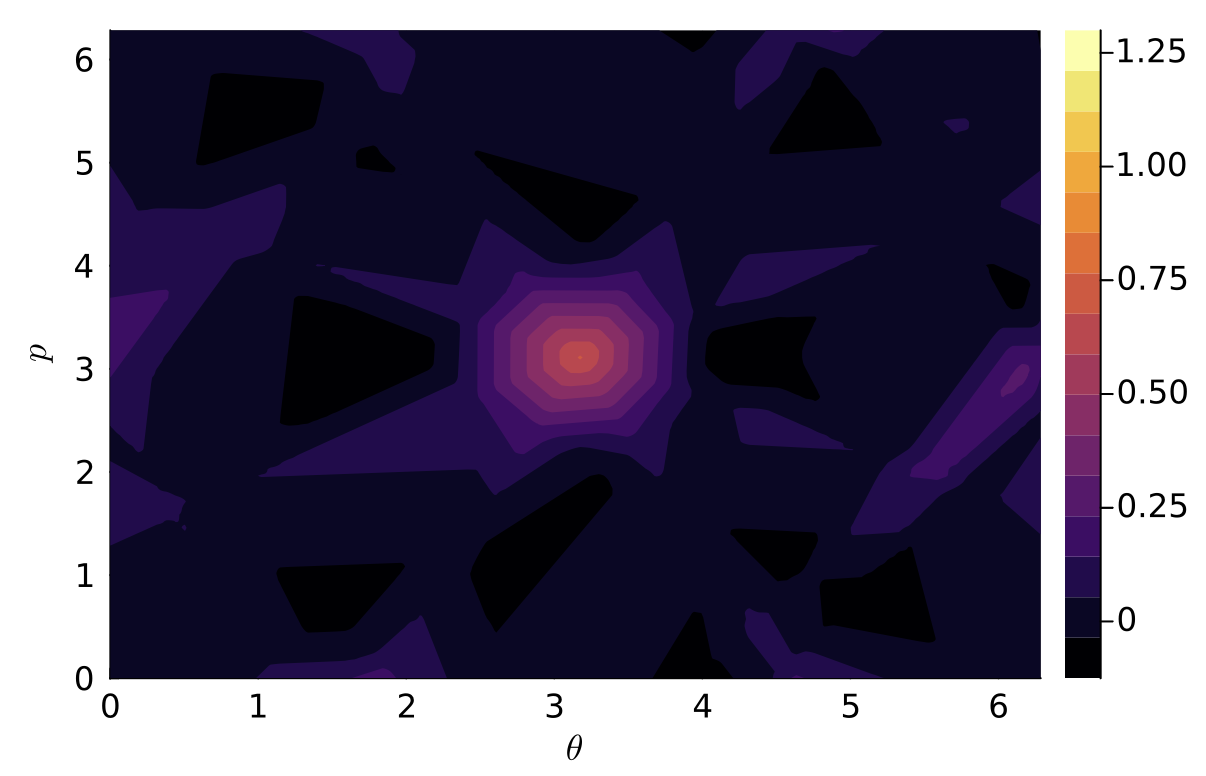}}\hfill 
    \subfloat[Two-layer neural network]{\includegraphics[width=6cm]{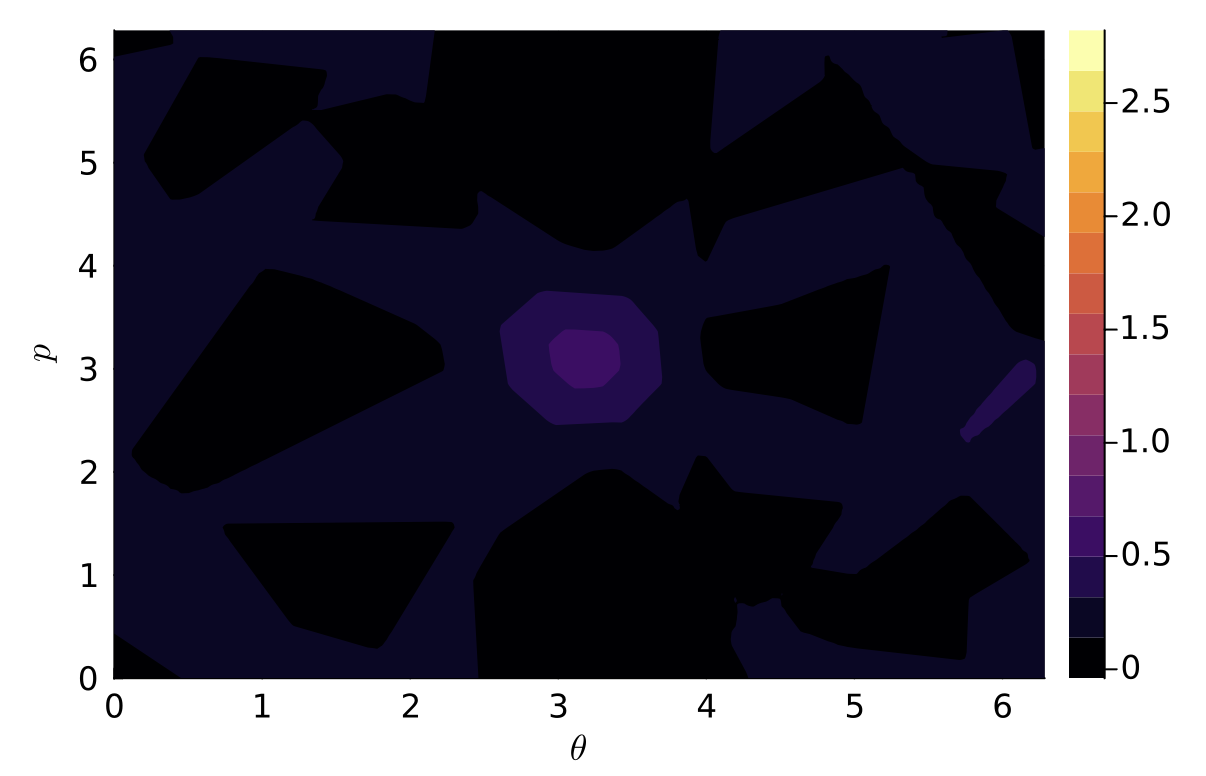}}\\
    \subfloat[Three-layer neural network]{\includegraphics[width=6cm]{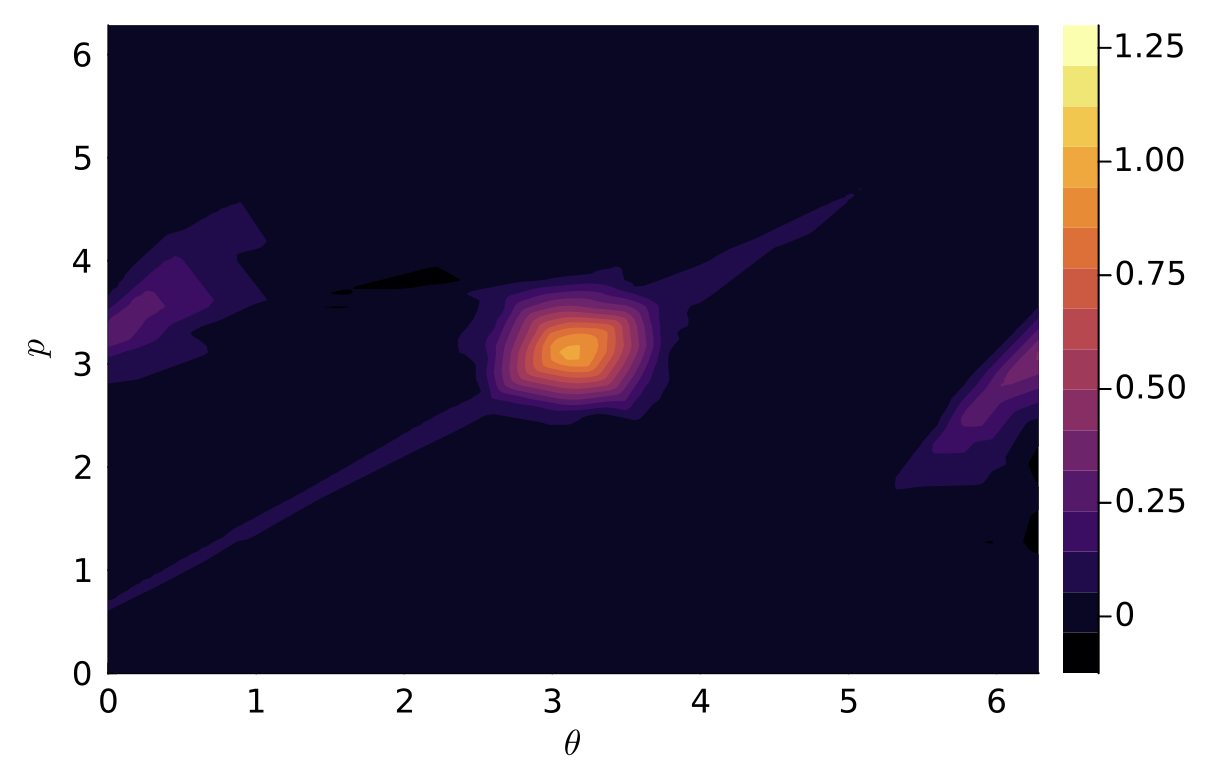}}\hfill 
    \subfloat[Three-layer neural network]{\includegraphics[width=6cm]{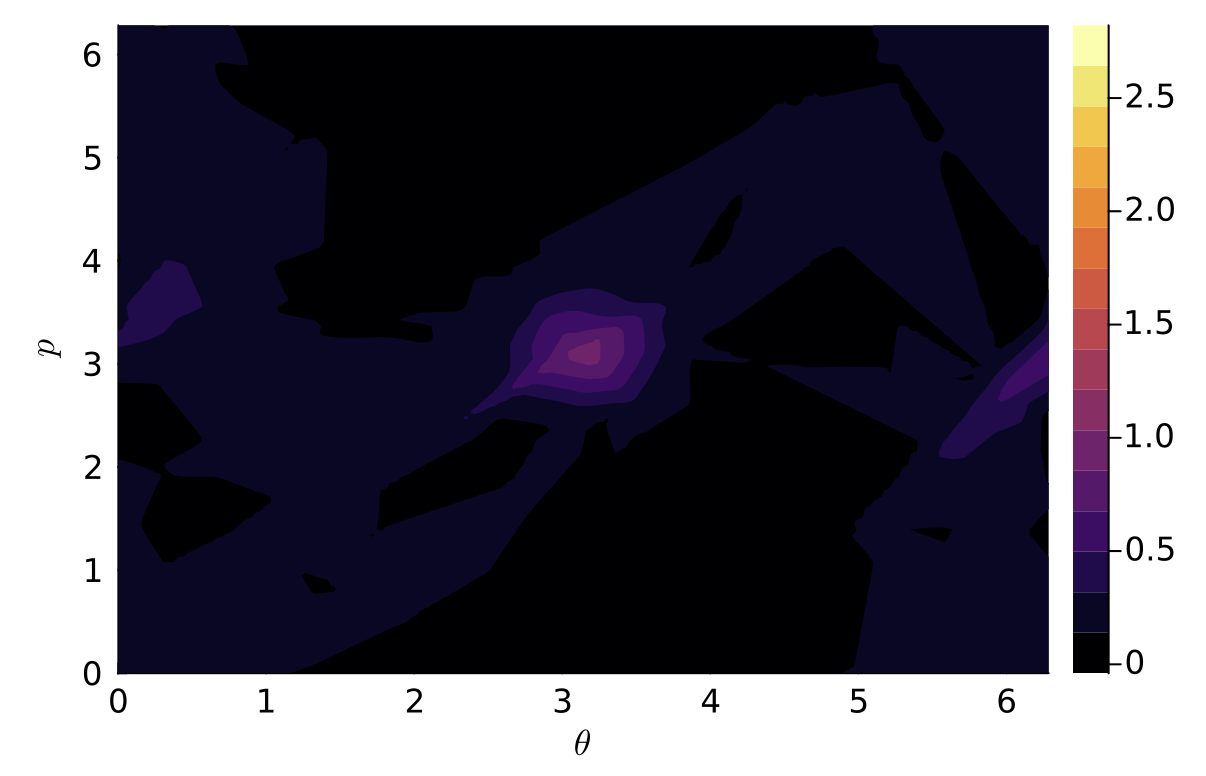}}
    \caption{Solutions for the standard map example with $\alpha = 0.5$ (left) and $\alpha = 0.9$ (right) computed by the truncated sum with $1{,}000$ terms and PINNs approximations.}
    \label{Figure: ANN approximation for the standard map}
\end{figure}

We also apply PINNs and RVPINNs to the standard map example, where the initial density is given by $f_0(\theta,p)=\cos^2(2\theta)\cos^2(2p)$ for $3\pi/4 < \theta < 5\pi/4$ and $3\pi/4 < p < 5\pi/4$, and $0$ elsewhere. Because $f_0$ is bounded with compact support, we have $f_0 \in L^2(\Omega)$. As $\alpha$ increases, Neumann series for the standard map example become more complex, as shown in Figure \ref{Figure: VPINNs steady state standard map}. To make neural networks converge to complex solutions faster, we employ a continuation technique: training the neural network to approximate solutions for small $\alpha$, then using the resulting model as an initialization for training with larger $\alpha$.
We train neural networks following the same approach described above. For the three-layer architecture, each hidden layer contains $32$ neurons with the ReLU activation function. Figure \ref{Figure: ANN approximation for the standard map} shows the approximate solutions obtained with PINNs for different values of $\alpha$. Similar results were observed for RVPINNs, which are therefore omitted for brevity.
These results demonstrate that PINNs and RVPINNs can capture the key features of solutions. Moreover, deeper neural networks yield better results due to their greater expressivity.

\begin{figure}[t]
    \centering

\tikzset{every picture/.style={line width=0.75pt}} %set default line width to 0.75pt        

\begin{tikzpicture}[x=0.75pt,y=0.75pt,yscale=-1,xscale=1]

%Image [id:dp7858426172005182] 
\draw (364.75,148.1) node  {\includegraphics[width=296.63pt,height=211.35pt]{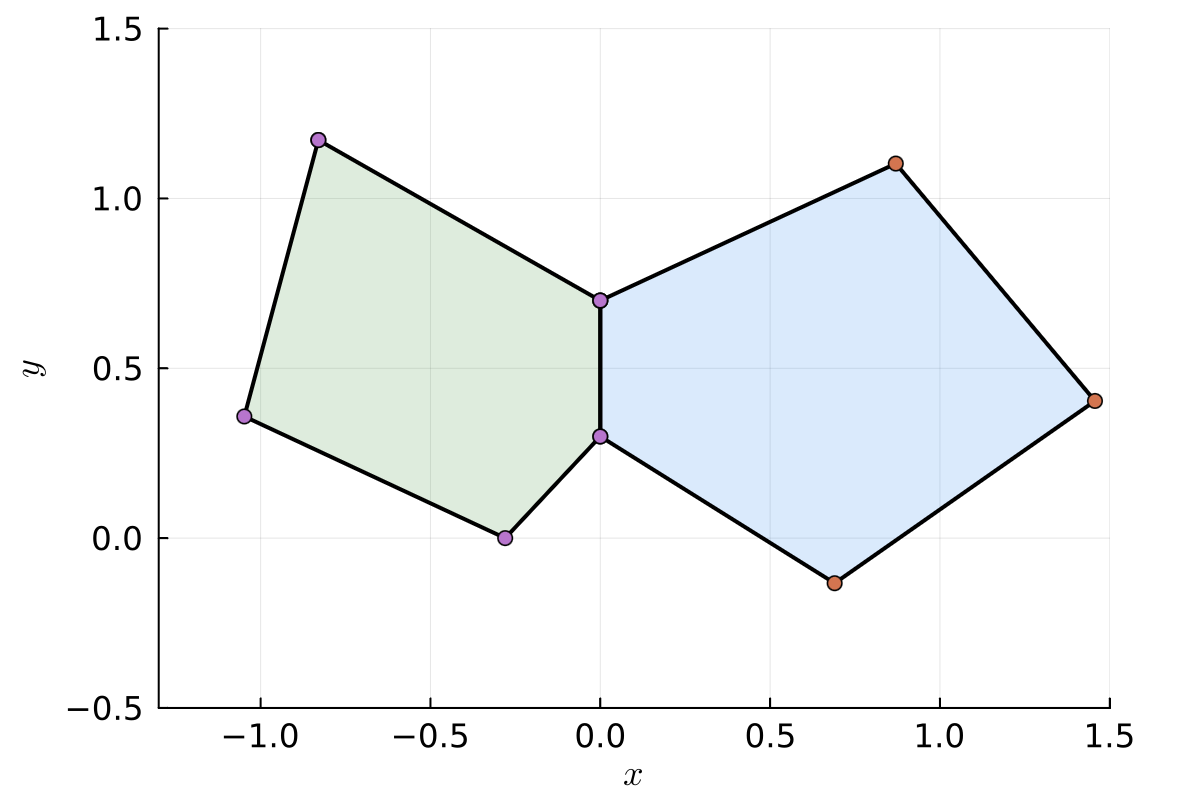}};
%Straight Lines [id:da2949227476500599] 
\draw  [dash pattern={on 4.5pt off 4.5pt}]  (373,118) -- (373,158) ;
%Straight Lines [id:da053011537383693685] 
\draw  [dash pattern={on 4.5pt off 4.5pt}]  (276,67) -- (254.5,150) ;
%Straight Lines [id:da3192441973396014] 
\draw  [dash pattern={on 4.5pt off 4.5pt}]  (254.5,150) -- (301.22,174.08) ;
\draw [shift={(303,175)}, rotate = 207.27] [color={rgb, 255:red, 0; green, 0; blue, 0 }  ][line width=0.75]    (10.93,-3.29) .. controls (6.95,-1.4) and (3.31,-0.3) .. (0,0) .. controls (3.31,0.3) and (6.95,1.4) .. (10.93,3.29)   ;
%Straight Lines [id:da7083394784574093] 
\draw  [dash pattern={on 4.5pt off 4.5pt}]  (373,158) -- (442,204) ;
%Straight Lines [id:da3872795579977384] 
\draw    (485,181) -- (494.37,152.9) ;
\draw [shift={(495,151)}, rotate = 108.43] [color={rgb, 255:red, 0; green, 0; blue, 0 }  ][line width=0.75]    (10.93,-3.29) .. controls (6.95,-1.4) and (3.31,-0.3) .. (0,0) .. controls (3.31,0.3) and (6.95,1.4) .. (10.93,3.29)   ;
%Straight Lines [id:da13365247669984104] 
\draw  [dash pattern={on 4.5pt off 4.5pt}]  (442,204) -- (478.9,176.2) ;
\draw [shift={(480.5,175)}, rotate = 143.01] [color={rgb, 255:red, 0; green, 0; blue, 0 }  ][line width=0.75]    (10.93,-3.29) .. controls (6.95,-1.4) and (3.31,-0.3) .. (0,0) .. controls (3.31,0.3) and (6.95,1.4) .. (10.93,3.29)   ;
%Straight Lines [id:da07163368812032334] 
\draw    (491,185) -- (510.9,170.19) ;
\draw [shift={(512.5,169)}, rotate = 143.34] [color={rgb, 255:red, 0; green, 0; blue, 0 }  ][line width=0.75]    (10.93,-3.29) .. controls (6.95,-1.4) and (3.31,-0.3) .. (0,0) .. controls (3.31,0.3) and (6.95,1.4) .. (10.93,3.29)   ;
%Straight Lines [id:da9386980966733383] 
\draw    (302,181) -- (331.32,162.08) ;
\draw [shift={(333,161)}, rotate = 147.17] [color={rgb, 255:red, 0; green, 0; blue, 0 }  ][line width=0.75]    (10.93,-3.29) .. controls (6.95,-1.4) and (3.31,-0.3) .. (0,0) .. controls (3.31,0.3) and (6.95,1.4) .. (10.93,3.29)   ;

% Text Node
\draw (407,160.4) node [anchor=north west][inner sep=0.75pt]    {$s$};
% Text Node
\draw (504,174.4) node [anchor=north west][inner sep=0.75pt]    {$p$};
% Text Node
\draw (267.25,111.9) node [anchor=north west][inner sep=0.75pt]    {$s'$};
% Text Node
\draw (357,90.4) node [anchor=north west][inner sep=0.75pt]    {$A$};
% Text Node
\draw (265,35.4) node [anchor=north west][inner sep=0.75pt]    {$B$};
% Text Node
\draw (411,207.4) node [anchor=north west][inner sep=0.75pt]    {$\Gamma _{1}$};
% Text Node
\draw (248,169.4) node [anchor=north west][inner sep=0.75pt]    {$\Gamma _{2}$};
% Text Node
\draw (434,122.4) node [anchor=north west][inner sep=0.75pt]    {$X_{1}$};
% Text Node
\draw (306,118.4) node [anchor=north west][inner sep=0.75pt]    {$X_{2}$};

\end{tikzpicture}

\caption{Two-cavity system: $\Gamma_1$ and $\Gamma_2$ are the boundaries of pentagons $X_1$ and $X_2$, respectively. When a ray lies in $X_1$, its arclength $s$ is measured anticlockwise from point $A$ along the boundary $\Gamma_1$. When a ray lies in $X_2$, its arclength $s$ is defined as $s = s'+|\Gamma_1|$, where $s'$ is the anticlockwise distance from point $B$ along the boundary $\Gamma_2$.}
\label{fig:two-cavity system}
\end{figure}

\begin{figure}[t]
    \centering

\tikzset{every picture/.style={line width=0.75pt}} %set default line width to 0.75pt        

\begin{tikzpicture}[x=0.75pt,y=0.75pt,yscale=-1,xscale=1]

%Image [id:dp7858426172005182] 
\draw (364.75,148.1) node  {\includegraphics[width=296.63pt,height=211.35pt]{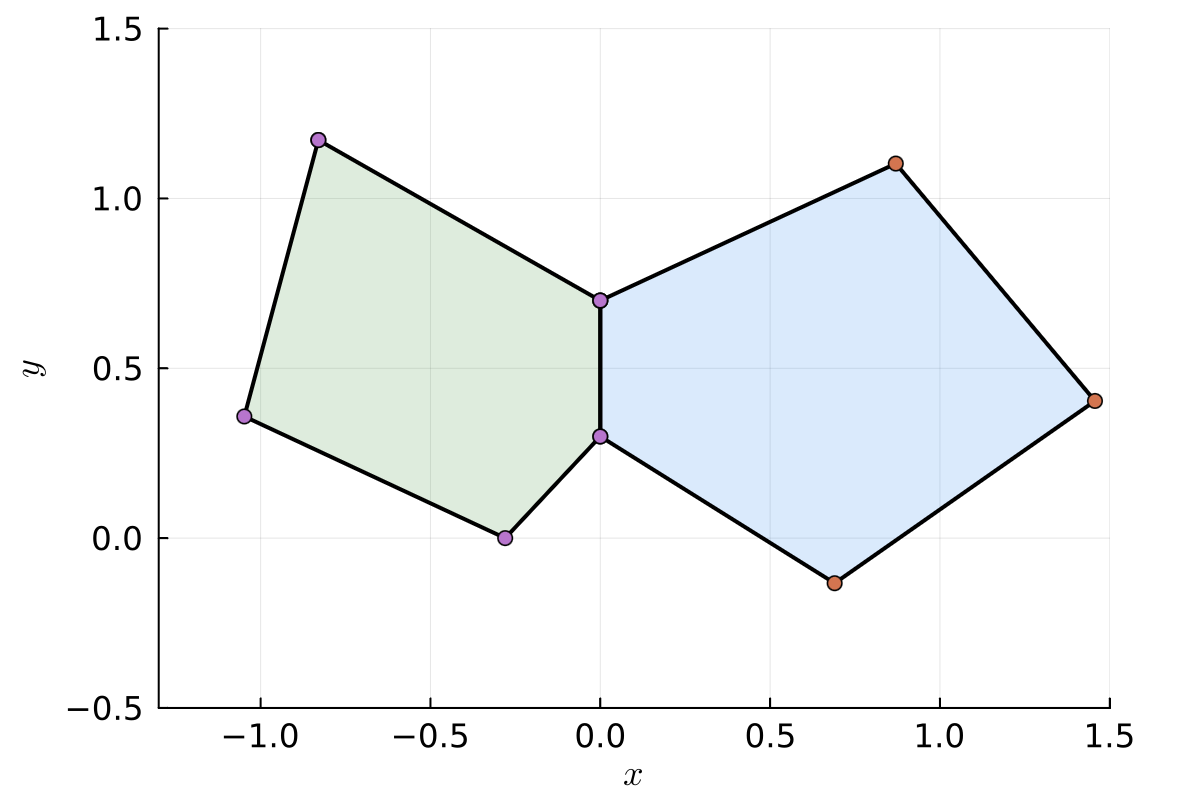}};
%Straight Lines [id:da21201277069569058] 
\draw [color={rgb, 255:red, 0; green, 0; blue, 0 }  ,draw opacity=1 ]   (451,205) -- (365.65,146.14) ;
\draw [shift={(364,145)}, rotate = 34.59] [color={rgb, 255:red, 0; green, 0; blue, 0 }  ,draw opacity=1 ][line width=0.75]    (10.93,-3.29) .. controls (6.95,-1.4) and (3.31,-0.3) .. (0,0) .. controls (3.31,0.3) and (6.95,1.4) .. (10.93,3.29)   ;
%Straight Lines [id:da2787195942180566] 
\draw [color={rgb, 255:red, 0; green, 0; blue, 0 }  ,draw opacity=1 ]   (364,145) -- (266.66,80.11) ;
\draw [shift={(265,79)}, rotate = 33.69] [color={rgb, 255:red, 0; green, 0; blue, 0 }  ,draw opacity=1 ][line width=0.75]    (10.93,-3.29) .. controls (6.95,-1.4) and (3.31,-0.3) .. (0,0) .. controls (3.31,0.3) and (6.95,1.4) .. (10.93,3.29)   ;
%Straight Lines [id:da10276633533986435] 
\draw [color={rgb, 255:red, 0; green, 0; blue, 0 }  ,draw opacity=1 ]   (265,79) -- (316.01,83.81) ;
\draw [shift={(318,84)}, rotate = 185.39] [color={rgb, 255:red, 0; green, 0; blue, 0 }  ,draw opacity=1 ][line width=0.75]    (10.93,-3.29) .. controls (6.95,-1.4) and (3.31,-0.3) .. (0,0) .. controls (3.31,0.3) and (6.95,1.4) .. (10.93,3.29)   ;
%Straight Lines [id:da602776732337887] 
\draw [color={rgb, 255:red, 0; green, 0; blue, 0 }  ,draw opacity=1 ]   (318,84) -- (356.17,168.18) ;
\draw [shift={(357,170)}, rotate = 245.61] [color={rgb, 255:red, 0; green, 0; blue, 0 }  ,draw opacity=1 ][line width=0.75]    (10.93,-3.29) .. controls (6.95,-1.4) and (3.31,-0.3) .. (0,0) .. controls (3.31,0.3) and (6.95,1.4) .. (10.93,3.29)   ;
%Straight Lines [id:da14286600870251787] 
\draw [color={rgb, 255:red, 0; green, 0; blue, 0 }  ,draw opacity=1 ]   (357,170) -- (259.75,115.97) ;
\draw [shift={(258,115)}, rotate = 29.05] [color={rgb, 255:red, 0; green, 0; blue, 0 }  ,draw opacity=1 ][line width=0.75]    (10.93,-3.29) .. controls (6.95,-1.4) and (3.31,-0.3) .. (0,0) .. controls (3.31,0.3) and (6.95,1.4) .. (10.93,3.29)   ;
%Straight Lines [id:da2859414095106477] 
\draw [color={rgb, 255:red, 0; green, 0; blue, 0 }  ,draw opacity=1 ]   (258,115) -- (363,121.87) ;
\draw [shift={(365,122)}, rotate = 183.74] [color={rgb, 255:red, 0; green, 0; blue, 0 }  ,draw opacity=1 ][line width=0.75]    (10.93,-3.29) .. controls (6.95,-1.4) and (3.31,-0.3) .. (0,0) .. controls (3.31,0.3) and (6.95,1.4) .. (10.93,3.29)   ;
%Straight Lines [id:da8462991822432503] 
\draw [color={rgb, 255:red, 0; green, 0; blue, 0 }  ,draw opacity=1 ]   (365,122) -- (511,129.89) ;
\draw [shift={(513,130)}, rotate = 183.09] [color={rgb, 255:red, 0; green, 0; blue, 0 }  ,draw opacity=1 ][line width=0.75]    (10.93,-3.29) .. controls (6.95,-1.4) and (3.31,-0.3) .. (0,0) .. controls (3.31,0.3) and (6.95,1.4) .. (10.93,3.29)   ;
%Straight Lines [id:da4691949824670192] 
\draw [color={rgb, 255:red, 0; green, 0; blue, 0 }  ,draw opacity=1 ]   (513,130) -- (512.07,158) ;
\draw [shift={(512,160)}, rotate = 271.91] [color={rgb, 255:red, 0; green, 0; blue, 0 }  ,draw opacity=1 ][line width=0.75]    (10.93,-3.29) .. controls (6.95,-1.4) and (3.31,-0.3) .. (0,0) .. controls (3.31,0.3) and (6.95,1.4) .. (10.93,3.29)   ;

% Text Node
\draw (411,207.4) node [anchor=north west][inner sep=0.75pt]    {$\Gamma _{1}$};
% Text Node
\draw (248,169.4) node [anchor=north west][inner sep=0.75pt]    {$\Gamma _{2}$};
% Text Node
\draw (416,162.4) node [anchor=north west][inner sep=0.75pt]    {$l_{1}$};
% Text Node
\draw (274,90.4) node [anchor=north west][inner sep=0.75pt]    {$l_{2}$};
% Text Node
\draw (276,63.4) node [anchor=north west][inner sep=0.75pt]    {$l_{3}$};
% Text Node
\draw (314,92.4) node [anchor=north west][inner sep=0.75pt]    {$l_{4}$};
% Text Node
\draw (302,144.4) node [anchor=north west][inner sep=0.75pt]    {$l_{5}$};
% Text Node
\draw (303,119.4) node [anchor=north west][inner sep=0.75pt]    {$l_{6}$};
% Text Node
\draw (428,107.4) node [anchor=north west][inner sep=0.75pt]    {$l_{7}$};
% Text Node
\draw (496,135.4) node [anchor=north west][inner sep=0.75pt]    {$l_{8}$};

\end{tikzpicture}

\caption{An example of a ray trajectory in the two-cavity system. Its trajectory length is $l = l_1+l_2+\ldots+l_8$.}
\label{fig:two-cavity trajectory}
\end{figure}

\subsection{Application: Interior densities in a two-cavity system}

Our proposed methods generate approximate functions for stationary boundary densities when $S$ is a boundary map. These boundary densities can then be used to compute the corresponding interior densities. In this section, we demonstrate applications of PINNs for approximating a stationary density within a complex domain. Specifically, we consider a two-cavity system as considered in \cite{bajars2017improved,chappell2011dynamical}. 

The phase-space coordinates on the boundary are denoted by $(s,p)$, where $s$ represents the arclength along the boundary and $p$ is the tangential component of the momentum vector (See Figure \ref{fig:two-cavity system}). Let $\Omega = (\Gamma_1 \cup \Gamma_2) \times(-1,1)$ be the phase space, and define the associated boundary map as a billiard map $S:\Omega \to \Omega$. It is well known that the billiard map is invertible almost everywhere, except at corner points where the reflection is not defined. Moreover, it is a Hamiltonian system \cite{cvitanovic2005chaos}, which implies that it preserves phase-space volume; equivalently, its Jacobian determinant is equal to one. Using Remark \ref{Remark: S diffeomorphism}, the Perron-Frobenius operator can be written as $\mathcal{P}f = f(S^{-1})$. 

Assume that the ray source is located on the third boundary segment of $\Gamma_1$. The corresponding initial boundary density is defined as $f_0(s,p) = (1-p^2)\sin(\pi(s-1.214)/0.936)$ for $1.214 < s < 2.150$ and $-1 < p < 1$, and $0$ elsewhere. The damping along each trajectory is modeled by the function $\alpha(l) = \exp(-0.5l)$, where $l$ denotes the trajectory length (see Figure \ref{fig:two-cavity trajectory}). 

Once the stationary boundary density $u(s,p)$ is obtained, the corresponding interior density can be computed by projecting $u$ onto the interior domain $X_i$ via the mapping
\begin{equation}\label{Eq: interior density}
    u_{X_i}(r) = \int_{\Gamma_i}u(s,p_s)\dfrac{\cos(\nu(r_s,r))}{|r-r_s|}ds,
\end{equation}
where $r \in X_i$ denotes a point inside the cavity (in Cartesian coordinates), $r_s$ is the Cartesian coordinate corresponding to the boundary point $s \in \Gamma_i$, and $\nu(r_s,r)$ is the angle between the inward-pointing normal vector at $s$ and the direction vector from $r_s$ to $r$. 

\begin{figure}[h]
    \subfloat[Truncated sum]{\includegraphics[width=6.5cm]{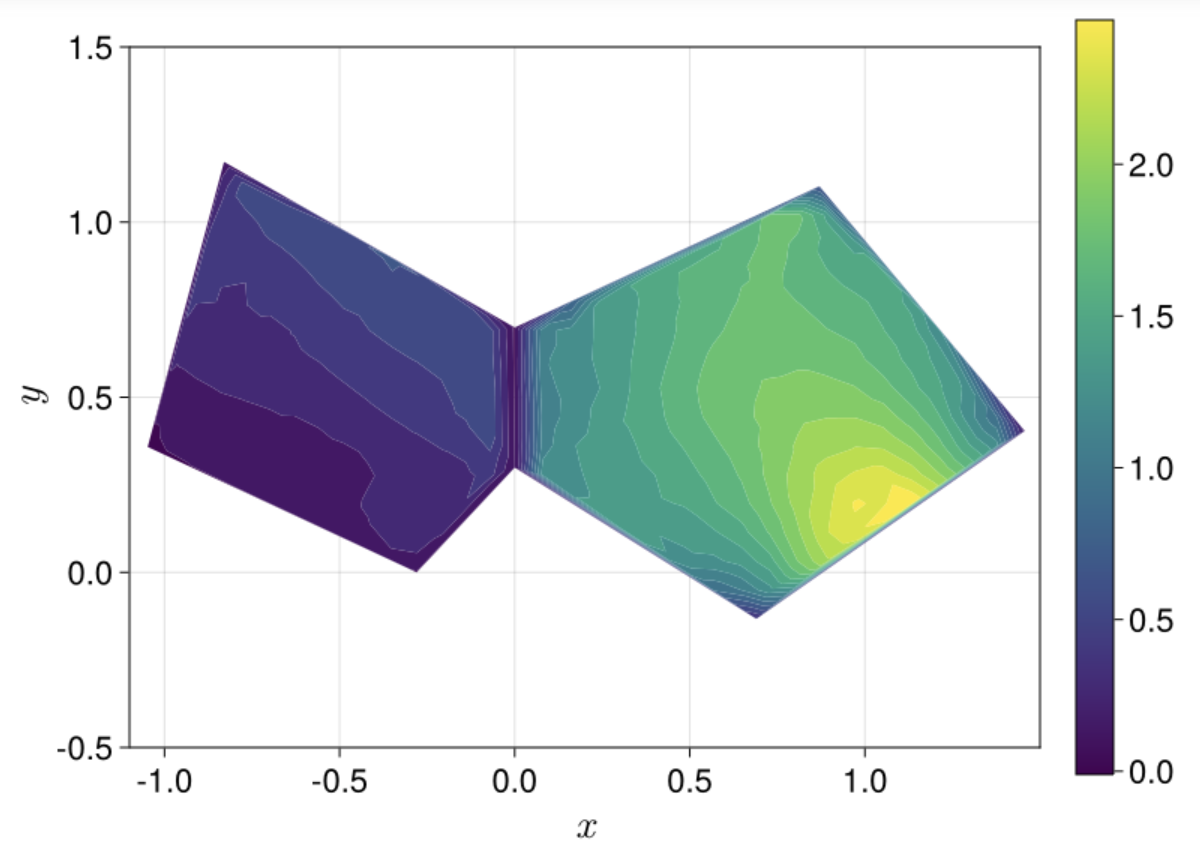}}\hfill 
    \subfloat[Fixed-grid-based method]{\includegraphics[width=6.5cm]{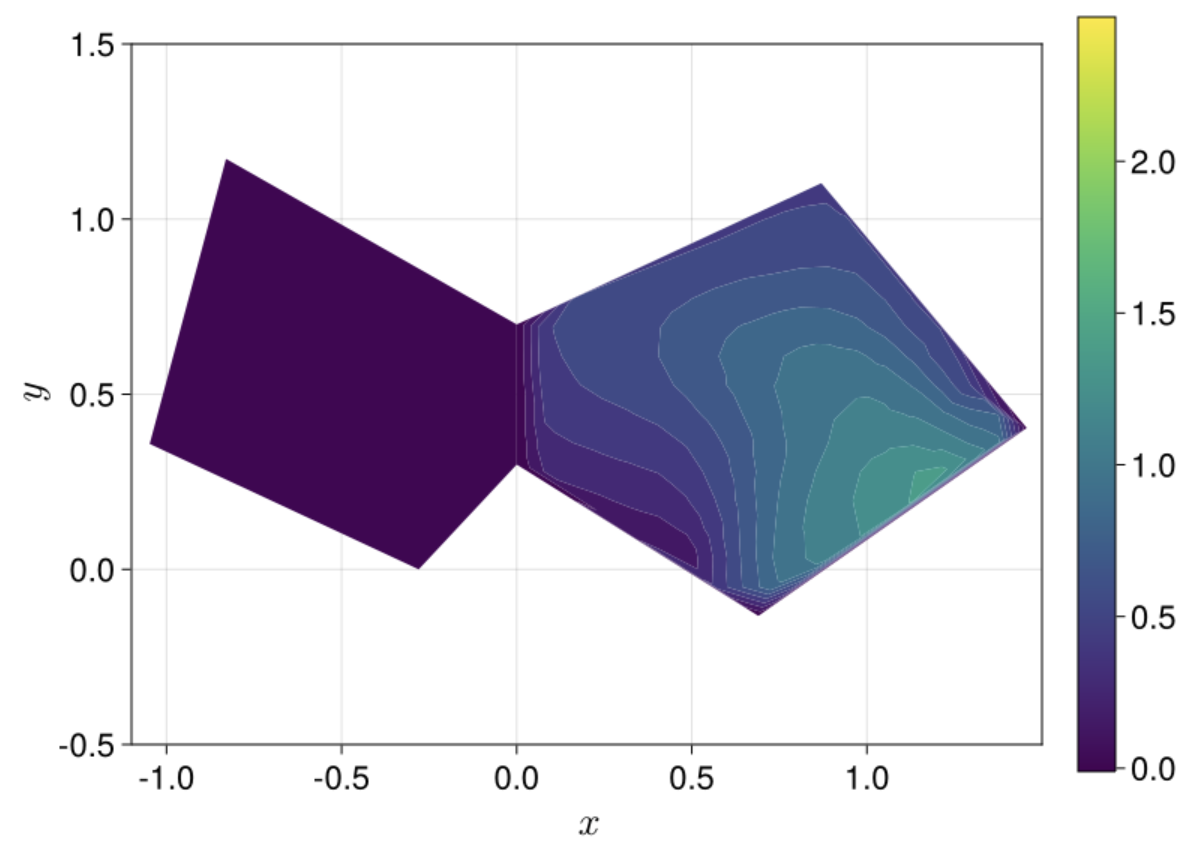}}\\
    \subfloat[PINNs with two-layer neural network]{\includegraphics[width=6.5cm]{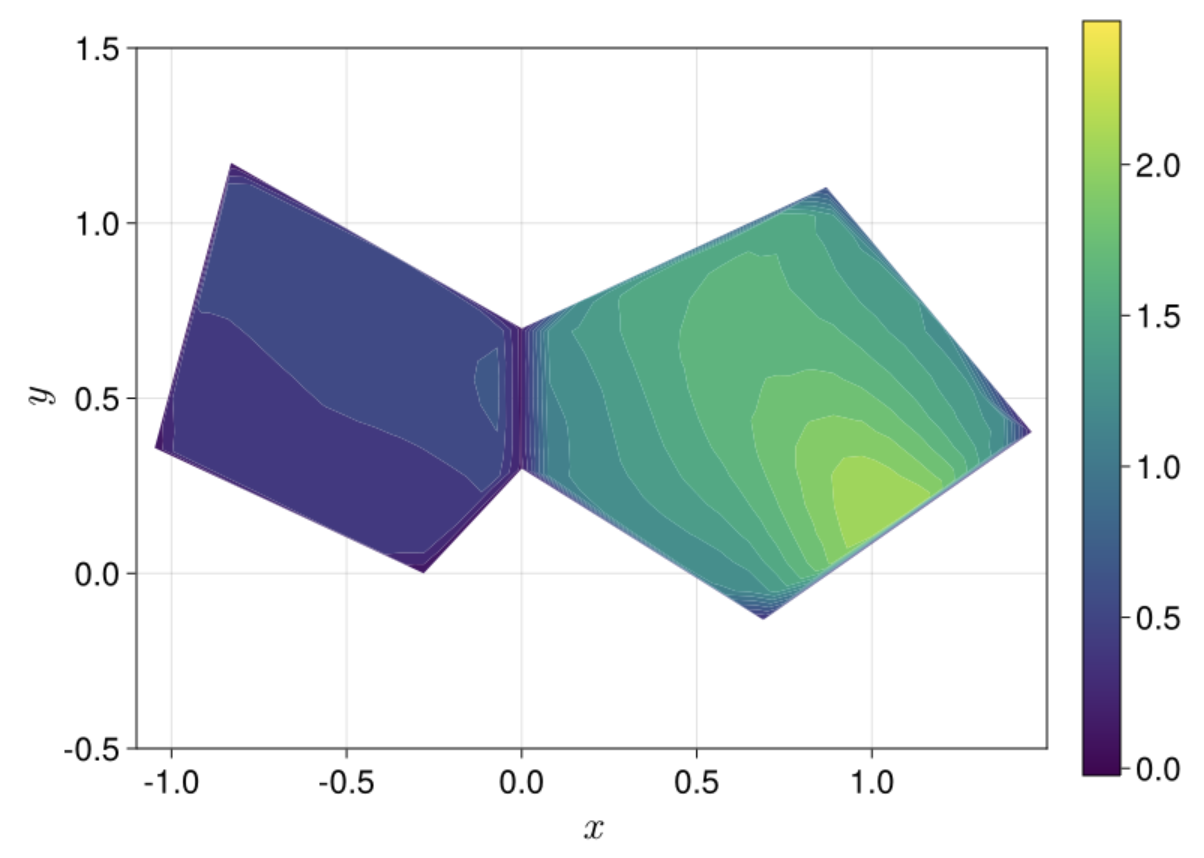}}\hfill 
    \subfloat[PINNs with three-layer neural network]{\includegraphics[width=6.5cm]{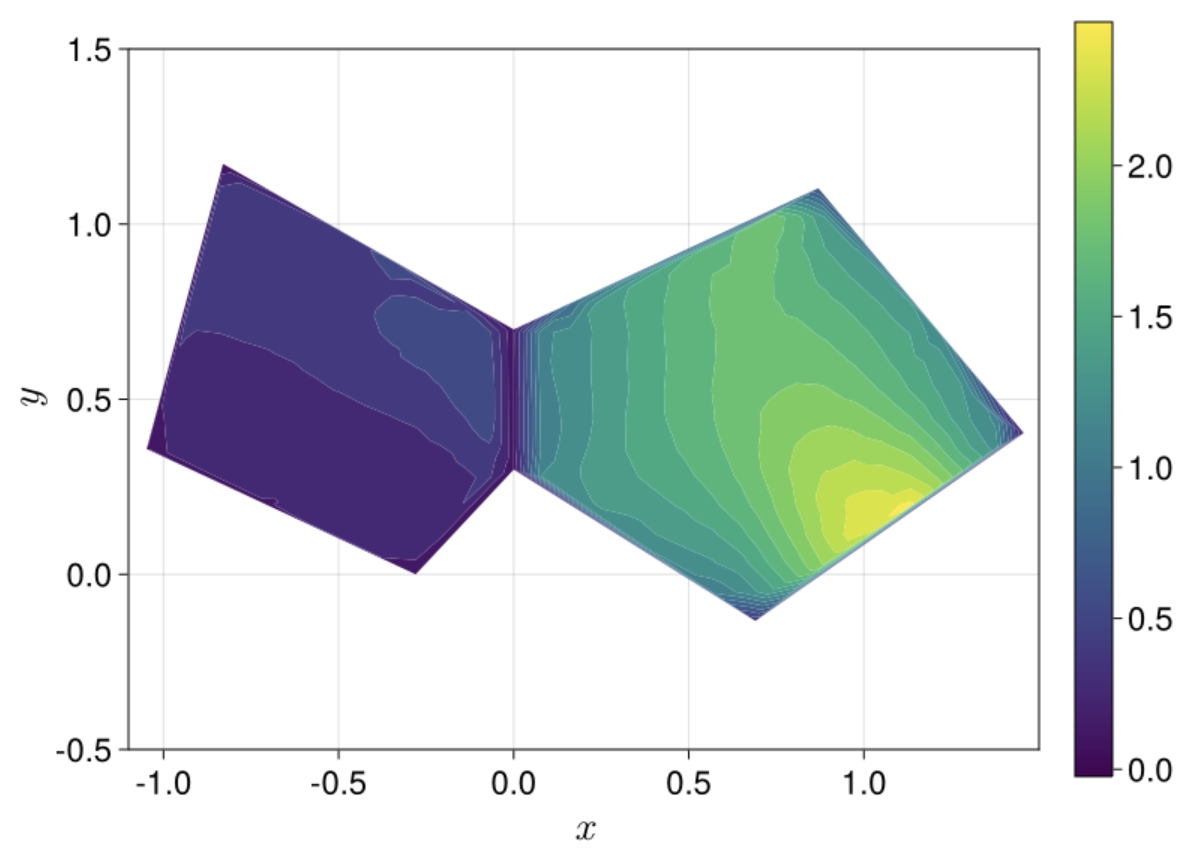}}
    \caption{Approximate stationary interior densities of the two-cavity system.}
    \label{Figure: two-cavity system}
\end{figure}

For the numerical results, Figure \ref{Figure: two-cavity system} shows the approximate stationary interior densities computed using (\ref{Eq: interior density}), where $u$ denotes the stationary boundary densities obtained from the truncated sum (``exact" solution), PINNs with two-layer and three-layer architectures, and the fixed-grid-based method. For PINNs, each hidden layer of the neural networks contains $64$ neurons with the ReLU activation function. For the fixed-grid-based method, we use piecewise linear spline functions defined on an $8\times 8$ uniform grid as basis functions. This choice is intentional to ensure a fair comparison with the two-layer (one-hidden-layer) neural network in terms of representational complexity: both approaches generate continuous, piecewise linear approximations and use $64$ basis functions. 

As shown in Figure \ref{Figure: two-cavity system}, the fixed-grid-based method cannot give fine details in the density distribution within the left pentagon, and the densities in the right pentagon exhibit noticeable deviations from the solution. In contrast, the two-layer neural network captures the main feature of the solution compared to the truncated sum. The three-layer neural network provides a more accurate approximation than the two-layer architecture. Overall, PINNs outperform the fixed-grid-based method.

\section{Conclusions}\label{sec: conclusions}

In this article, we investigated Neumann series problems associated with non-expansive Perron-Frobenius operators under the $L^p$-norm, with a focus on the $L^2$-setting. We established the well-posedness of the corresponding variational formulations in $L^2$ by showing that they satisfy the conditions of the Lax-Milgram theorem, thereby ensuring the existence and uniqueness of the solution. We proposed two neural network methods for approximating solutions of Neumann series problems. PINNs were employed to approximate solutions in their strong form, while RVPINNs were used to approximate solutions in the corresponding variational problems.  Additionally, we provided error estimates for both methods in terms of quasi-minimizers. Numerical experiments demonstrated that both methods are effective for solving Neumann series problems and outperform the standard fixed-grid-based method. 

Further extensions of this work include applying our methods to Neumann series approximations in higher-dimensional settings. In addition, the theoretical results in this paper, such as a priori error estimates based on the local Fortin's condition, could be applied to establish reliable error bounds for solving PDEs with neural networks.

\section*{Acknowledgments}
T. Udomworarat was supported by the Royal Thai Government Scholarship under the Office of Educational Affairs, the Royal Thai Embassy in London. The research by I. Brevis and K. G. van der Zee was supported by the Engineering and Physical Sciences Research Council (EPSRC), UK, under Grant EP/W010011/1. S. Rojas's work was supported by the Chilean grant ANID FONDECYT No. 1240643, and by the National Center for Artificial Intelligence CENIA FB210017, Basal ANID.

\appendix
\section{Well-posedness of variational problem} \label{section: well-posedness}
In the following, we show that the variational problem \cref{steady state variational problem} satisfies the assumptions of Banach-Ne$\check{\textrm{c}}$as-Bab$\check{\textrm{u}}$ska Theorem, so the existence and uniqueness of the solution are verified. 

\begin{thm}[Well-posedness of variational problem in $L^p-L^q$ spaces]\label{thm:Well-posedness Lp-Lq}
  If $U=L^p(\Omega)$ and $V=L^q(\Omega)$, where $1<p,q<\infty$ and $\dfrac{1}{p}+\dfrac{1}{q}=1$. The variational problem \cref{steady state variational problem}  satisfies the following:
  \begin{enumerate}
      \item Boundedness of $b$: $|b(u,v)| \leq (1+\alpha)\|u\|_{L^p}\|v\|_{L^q}, \quad \forall u \in L^p(\Omega),v \in L^q(\Omega)$,
      \item Boundedness of $l$: $|l(v)| \leq \|f_0\|_{L^p}\|v\|_{L^q}, \quad \forall v \in L^q(\Omega)$,
      \item Inf-sup stability: $\sup\limits_{0\neq v \in L^q(\Omega)}\frac{b(u,v)}{\|v\|_{L^q}} \geq (1-\alpha) \|u\|_{L^p}, \quad \forall u \in L^p(\Omega)$,
      \item Adjoint injectivity: $b(u,v) = 0 \quad \forall u \in L^p(\Omega) \Longrightarrow v = 0$.
  \end{enumerate}
  Hence, the variational problem \cref{steady state variational problem} has a unique solution.
\end{thm}

\begin{pf}
    1. Boundedness of $b$: For $u \in L^p(\Omega)$ and $v \in L^q(\Omega)$, by using H{\"o}lder's inequality, the triangle inequality, and the assumption (A2), we have
    \begin{align*}
          |b(u,v)| &= \left|\int_{\Omega}(u-\alpha\mathcal{P}u) v d\mu \right| \leq \int_{\Omega}|(u-\alpha\mathcal{P}u) v| d\mu \leq \|u-\alpha\mathcal{P}u\|_{L^p}\|v\|_{L^q} \\
          &\leq \left(\|u\|_{L^p} + \alpha\|\mathcal{P}u\|_{L^p} \right)\|v\|_{L^q} \leq (1+\alpha)\|u\|_{L^p}\|v\|_{L^q}.
      \end{align*}
      
    2. Boundedness of $l$: Since $f_0 \in L^p(\Omega),v \in L^q(\Omega)$, by using H{\"o}lder's inequality, we have
    \begin{align*}
        |l(v)| = \left|\int_{\Omega} f_0 v d\mu\right| \leq \|f_0\|_{L^p} \|v\|_{L^q}.
    \end{align*}
    
    3. Inf-sup stability: For $u \in L^p(\Omega)$, by substituting $v=(u-\alpha\mathcal{P}u)^{p-1} \in L^q(\Omega)$, and using the reverse triangle inequality and assumption (A2), we have
    \begin{align*}
        \sup_{0\neq v \in L^q(\Omega)}\frac{b(u,v)}{\|v\|_{L^q}} 
        &= \sup_{0\neq v \in L^q(\Omega)}\frac{\int_{\Omega}(u-\alpha\mathcal{P}u)vd\mu}{\|v\|_{L^q}} \\
        &\geq \frac{\int_{\Omega}(u-\alpha\mathcal{P}u)^pd\mu}{\|(u-\alpha\mathcal{P}u)^{p-1}\|_{L^q}} 
        = \frac{\|u-\alpha\mathcal{P}u\|_{L^p}^p}{\|u-\alpha\mathcal{P}u\|_{L^p}^{p/q}} \\
        &= \|u-\alpha\mathcal{P}u\|_{L^p} 
        \geq \|u\|_{L^p} - \alpha\|\mathcal{P}u\|_{L^p} 
        \geq (1-\alpha)\|u\|_{L^p}.
    \end{align*}

    4. Adjoint injectivity: Let $v \in L^q(\Omega)$ and assume that $b(u,v) = 0 \quad \forall u \in L^p(\Omega)$. By choosing $u=v^{q-1} \in L^p(\Omega)$, and using H{\"o}lder's inequality and the assumption (A2), we have
    \begin{align*}
        0 &= \int_{\Omega}(v^{q-1}-\alpha\mathcal{P}v^{q-1})vd\mu 
        = \int_{\Omega} v^q d\mu - \alpha\int_{\Omega}(\mathcal{P}v^{q-1})vd\mu \\
        &\geq \|v\|_{L^q}^q - \alpha\|\mathcal{P}v^{q-1}\|_{L^p}\|v\|_{L^q}
        \geq \|v\|_{L^q}^q - \alpha\|v^{q-1}\|_{L^p}\|v\|_{L^q} \\ 
        &= \|v\|_{L^q}^q - \alpha\|v\|_{L^q}^{q/p}\|v\|_{L^q} 
        = (1-\alpha)\|v\|_{L^q}^q.
    \end{align*}
    Hence, $\|v\|_{L^q} = 0$, i.e., $v = 0$.
\end{pf}

\section{Equivalent form of RVPINNs loss function}\label{section: equivalent form of RVPINNs loss}

For $u_{\theta} \in \mathcal{M}_n \subset L^p(\Omega)$ and $v_M \in V_M \subset L^q(\Omega)$, we define 
$$r(u_{\theta},v_M):=l(v_M)-b(u_{\theta},v_M).$$ 
Suppose that $\tilde{g}_{\theta} \in V_M$ is the solution of the following problem:
\begin{equation}\label{Eq: duality map equation}
    \left<J_M(\tilde{g}_{\theta}),v_M \right>_{L^p,L^q} = r(u_{\theta},v_M), \quad \textrm{for all } v_M \in V_M,
\end{equation}
where $\left<\cdot,\cdot \right>_{L^p,L^q}$ denotes the duality pairing and $J_M$ is the duality map on a subspace $V_M$ defined by $J_M(g) = \|g\|_{L^q}^{2-q} |g|^{q-1}\mathrm{sign}(g)$ (see  \cite{muga2020discretization}).
The following expressions are equivalent:
\begin{enumerate}
    \item $\sup\limits_{0 \neq v_M \in V_M} \dfrac{l(v_M)-b(u_{\theta},v_M)}{\|v_M\|_{L^q}}$,
    \item $\sup\limits_{0 \neq v_M \in V_M} \dfrac{b(u - u_{\theta},v_M)}{\|v_M\|_{L^q}}$, where $u$ is the solution of problem \cref{steady state variational problem},
    \item $\|J_M(\tilde{g}_{\theta})\|_{L^p}$.
\end{enumerate}

Indeed, the equivalence of 1. and 2. follows from $b(u,v_M) = l(v_M)~ \forall v_M \in V_M$ and the linearity of $b$. The equivalence of 1. and 3. follows from
\begin{align*}
    \sup_{0 \neq v_M \in V_M} \frac{l(v_M)-b(u_{\theta},v_M)}{\|v_M\|_{L^q}}  &= \sup_{0 \neq v_M \in V_M} \frac{r(u_{\theta},v_M)}{\|v_M\|_{L^q}} \\
    &= \sup_{0 \neq v_M \in V_M} \frac{\left< J_M(\tilde{g}_{\theta}),v_M\right>_{L^p,L^q}}{\|v_M\|_{L^q}} = \|J_M(\tilde{g}_{\theta})\|_{L^p},
\end{align*}
where the second equality follows from \cref{Eq: duality map equation} and the third one from the fact that $v_M = J_M(\tilde{g}_{\theta})$ is the supremizer.

Thus, the RVPINNs loss function can be written in three different forms:
\begin{enumerate}
    \item $\mathcal{L}_{\mathrm{RVPINNs}}(u_{\theta}) = \sup\limits_{0 \neq v_M \in V_M} \dfrac{l(v_M)-b(u_{\theta},v_M)}{\|v_M\|_{L^q}}$,
    \item $\mathcal{L}_{\mathrm{RVPINNs}}(u_{\theta}) = \sup\limits_{0 \neq v_M \in V_M} \dfrac{b(u - u_{\theta},v_M)}{\|v_M\|_{L^q}}$,
    \item $\mathcal{L}_{\mathrm{RVPINNs}}(u_{\theta}) = \|J_M(\tilde{g}_{\theta})\|_{L^p}$.
\end{enumerate}

Now, consider the case where $u$ is the solution of the variational problem \cref{steady state variational problem} with $U=V=L^2(\Omega)$, and $u_{\theta} \in \mathcal{M}_n \subset L^2(\Omega)$.  Let $V_M$ be a finite-dimensional subspace of $L^2(\Omega)$ with an orthonormal basis $\{ g_1,\ldots,g_M \}$. In this setting, the duality paring in \cref{Eq: duality map equation} corresponds to the $L^2$-inner product, and the duality map $J_M(\tilde{g}_{\theta})$ is the Riesz representation in $V_M$ of $r(u_{\theta},\cdot)$.
Suppose that $g_{\theta} \in V_M$ is the solution of the following Galerkin problem:
\begin{equation*}
    \left<g_{\theta},v_M \right>_{L^2} = r(u_{\theta},v_M), \quad \textrm{for all } v_M \in V_M.
\end{equation*} 
Since $g_{\theta} \in V_M$, it can be written as
\begin{equation}\label{expansion of g}
    g_{\theta} = \sum_{m=1}^M \alpha_m g_m, \textrm{ where } \alpha_m \in \mathbb{R}, m = 1,2,\ldots, M.
\end{equation}
For each $m \in \{1,2,\ldots,M\}$, we compute the coefficients $\alpha_m$ by
\begin{equation*}
    \alpha_m =  \left<g_{\theta},g_m \right>_{L^2} = r(u_{\theta},g_m).
\end{equation*}
Substituting these into \cref{expansion of g}, we obtain
\begin{equation*}
    g_{\theta} = \sum_{m=1}^M r(u_{\theta},g_m) g_m.
\end{equation*}
Taking the $L^2$-norm on both sides, we have
\begin{equation*}
    \|g_{\theta}\|_{L^2} = \sqrt{\sum_{m=1}^M r^2(u_{\theta},g_m)}. 
\end{equation*}
Thus, the RVPINNs loss function for solutions in $L^2(\Omega)$ is 
\begin{equation*}
    \mathcal{L}_{\mathrm{RVPINNs}}(u_{\theta}) = \sqrt{\sum\limits_{m=1}^M \left(\int_{\Omega} (f_0 - u_{\theta} + \alpha\mathcal{P}u_{\theta})g_m d\mu\right)^2}.
\end{equation*}

\section{Error estimate proofs for PINNs and RVPINNs}\label{section: error estimates proof}

\underline{Proof of Theorem \ref{Thm: error estimate PINNs}}

For arbitrary $u_{\theta} \in \mathcal{M}_n$, by using the assumption (A2), the reverse triangle inequality, the definition of a quasi-minimizer, and the triangle inequality, we have
\begin{align*}
    (1-\alpha)\|u-u_{\theta^*}\|_{L^p} &\leq \|u-u_{\theta^*}\|_{L^p} - \alpha\|\mathcal{P}(u-u_{\theta^*})\|_{L^p}\\
    &\leq \|(u-u_{\theta^*}) - \alpha\mathcal{P}(u-u_{\theta^*})\|_{L^p} \\
    &= \|f_0 - u_{\theta^*} + \alpha\mathcal{P}u_{\theta^*}\|_{L^p} \\
    &\leq \|f_0 - u_{\theta} + \alpha\mathcal{P}u_{\theta}\|_{L^p} + \delta_n \\
    &= \|(u-u_{\theta}) - \alpha\mathcal{P}(u-u_{\theta})\|_{L^p} + \delta_n \\
    &\leq \|u-u_{\theta^*}\|_{L^p} + \alpha\|\mathcal{P}(u-u_{\theta^*})\|_{L^p} + \delta_n \\
    &\leq (1+\alpha)\|u-u_{\theta}\|_{L^p} + \delta_n.
\end{align*}
Dividing by $1-\alpha$ and taking the infimum over all possible $u_{\theta} \in \mathcal{M}_n$ yields the theorem.

In RVPINNs, the local Fortin's condition provides the following result, which is necessary for proving Theorem \ref{Thm: Error estimates variational problems}. 

\begin{lem}\label{Local Fortin's Lemma}
    Let $\delta_n > 0$ and $u_{\theta^*} \in \mathcal{M}_n$ be a quasi-minimizer of \cref{Eq: loss function VPINNs sup} satisfying \cref{Eq: Quasi-minimizer}. If the local Fortin's condition is satisfied, it holds:
    \begin{equation*}
        \|u_{\theta^*} - u_{\theta}\|_{L^p} \leq \frac{C_{\Pi}}{1-\alpha}\left(\sup_{0 \neq v_M \in V_M} \frac{b(u_{\theta^*} - u_{\theta},v_M)}{\|v_M\|_{L^q}}\right), \quad \forall u_{\theta} \in \mathcal{M}_n^{\theta^*,R}.
    \end{equation*}
\end{lem}
\begin{pf}
    By using the local Fortin's condition and the inf-sup stability in Theorem \ref{thm:Well-posedness Lp-Lq}, we have
    \begin{align*}
        \sup_{0\neq v_M \in V_M} \frac{b(u_{\theta^*} - u_{\theta},v_M)}{\|v_M\|_{L^q}} 
        &\geq \sup_{0\neq v \in L^q(\Omega)} \frac{b(u_{\theta^*} - u_{\theta},\Pi_{\theta}v)}{\|\Pi_{\theta}v\|_{L^q}} \\
        &= \sup_{0\neq v \in L^q(\Omega)} \frac{b(u_{\theta^*} - u_{\theta},v)}{\|\Pi_{\theta}v\|_{L^q}} \\
        &\geq \frac{1}{C_{\Pi}}\sup_{0\neq v \in L^q(\Omega)} \frac{b(u_{\theta^*} - u_{\theta},v)}{\|v\|_{L^q}} \\
        &\geq \frac{1-\alpha}{C_{\Pi}}\|u_{\theta^*} - u_{\theta}\|_{L^p}. 
    \end{align*}
    Multiplying both sides by $\dfrac{C_{\Pi}}{1-\alpha}$ yields the lemma.
\end{pf}

\underline{Proof of Theorem \ref{Thm: Error estimates variational problems}} 

For any $u_{\theta} \in \mathcal{M}_n^{\theta^*,R}$, by using the triangle inequality, Lemma \ref{Local Fortin's Lemma}, the definition of a quasi-minimizer, and the boundedness of $b(\cdot,\cdot)$, we have
\begin{align*}
    \|u-u_{\theta^*}\|_{L^p} &\leq \|u_{\theta^*} - u_{\theta}\|_{L^p} + \|u-u_{\theta}\|_{L^p} \\
    &\leq \dfrac{C_{\Pi}}{1-\alpha}\left( \sup_{0 \neq v_M \in V_M} \frac{b(u_{\theta^*} - u_{\theta},v_M)}{\|v_M\|_{L^q}} \right) + \|u-u_{\theta}\|_{L^p} \\
    &\leq \dfrac{C_{\Pi}}{1-\alpha}\left( \sup_{0 \neq v_M \in V_M} \frac{b(u - u_{\theta},v_M)}{\|v_M\|_{L^q}} + \sup_{0 \neq v_M \in V_M} \frac{b(u - u_{\theta^*},v_M)}{\|v_M\|_{L^q}} \right) \\ &~~~~+ \|u-u_{\theta}\|_{L^p} \\
    &\leq \dfrac{C_{\Pi}}{1-\alpha}\left( 2\sup_{0 \neq v_M \in V_M} \frac{b(u - u_{\theta},v_M)}{\|v_M\|_{L^q}} + \delta_n \right) + \|u-u_{\theta}\|_{L^p} \\
    &\leq \dfrac{2C_{\Pi}}{1-\alpha}(1+\alpha)\|u-u_{\theta}\|_{L^p} + \left(\dfrac{C_{\Pi}}{1-\alpha}\right)\delta_n + \|u-u_{\theta}\|_{L^p} \\
    &= \left( 1 + \dfrac{2C_{\Pi}}{1-\alpha}(1+\alpha) \right)\|u-u_{\theta}\|_{L^p} + \left( \dfrac{C_{\Pi}}{1-\alpha} \right)\delta_n.
\end{align*}
Taking the infimum over all possible $u_{\theta} \in \mathcal{M}_n^{\theta^*,R}$ yields the theorem.

\bibliographystyle{elsarticle-num} 
\bibliography{references}

\end{document}